\documentclass[12pt]{article}
\usepackage{amsmath,amsxtra,latexsym,amsthm,amssymb,amscd,amsfonts}
\usepackage[portrait, top=3.5cm, bottom=3cm, left=3.5cm, right=2cm] {geometry}
\usepackage{graphicx}
\usepackage{color}
\usepackage{enumerate}

\theoremstyle{plain}
\newtheorem{definition}{\bf Definition}[section]
\newtheorem{theorem}{\bf Theorem}[section]

\def\tto{\;{\lower 1pt \hbox{$\rightarrow$}}\kern -10pt
\hbox{\raise 2pt \hbox{$\rightarrow$}}\;}

\def\R{{\rm I\!R}}
\def\N{{\rm I\!N}}

\def\co{\mbox{\rm co}\,}

\def\gph{\mbox{\rm gph}\,}
\def\epi{\mbox{\rm epi}\,}

\begin{document}
\pagestyle{myheadings}

\newtheorem{Theorem}{Theorem}[section]
\newtheorem{Proposition}[Theorem]{Proposition}
\newtheorem{Remark}[Theorem]{Remark}
\newtheorem{Lemma}[Theorem]{Lemma}
\newtheorem{Corollary}[Theorem]{Corollary}
\newtheorem{Definition}[Theorem]{Definition}
\newtheorem{Example}[Theorem]{Example}
\renewcommand{\theequation}{\thesection.\arabic{equation}}
\normalsize

\setcounter{equation}{0}

\title{Analyzing a Maximum Principle for Finite Horizon State Constrained Problems via Parametric Examples. Part 2: Problems with Bilateral State Constraints\footnote{In this research,  Vu Thi Huong and Nguyen Dong Yen were supported by National Foundation for Science \& Technology Development (Vietnam)  under grant number 101.01-2018.308.}}

\author{V.T.~Huong\footnote{Institute of Mathematics, Vietnam Academy of Science and Technology, 18 Hoang
Quoc Viet, Hanoi 10307, Vietnam; email:
vthuong@math.ac.vn; huong263@gmail.com.},\ \, J.-C.~Yao\footnote{Research Center for Interneural Computing, China Medical University 
Hospital, China Medical University, Taichung 40402, Taiwan; Email: yaojc@mail.cmu.edu.tw.},\ \, and\ \, N.D. Yen\footnote{Institute
of Mathematics, Vietnam Academy of Science and Technology, 18 Hoang
Quoc Viet, Hanoi 10307, Vietnam; email: ndyen@math.ac.vn.}}

\maketitle
\date{}%\small\today}
\centerline{\textit{(Dedicated to Professor Gue Myung Lee on the occasion of his 65th birthday)}}

\bigskip
\begin{quote}
\noindent {\bf Abstract.}
In the present paper, the maximum principle for finite horizon state constrained problems from the book by R. Vinter [\textit{Optimal Control}, Birkh\"auser, Boston, 2000; Theorem~9.3.1] is analyzed via parametric examples. The latter has origin in a recent paper by  V.~Basco, P.~Cannarsa, and H.~Frankowska, and resembles the optimal growth problem in mathematical economics. The solution existence of these parametric examples is established by invoking Filippov's existence theorem for Mayer problems. Since the maximum principle is only a necessary condition for local optimal processes, a large amount of additional investigations is needed to obtain a comprehensive synthesis of finitely many processes suspected for being local minimizers. Our analysis not only helps to understand the principle in depth, but also serves as a sample of applying it to meaningful prototypes of  economic optimal growth models. Problems with unilateral state constraints have been studied in Part 1 of the paper. Problems with bilateral state constraints are addressed in this Part 2. 

\noindent {\bf Keywords:}\  Finite horizon optimal control problem, state constraint,  maximum principle, solution existence theorem, function of bounded variation, Borel measurable function, Lebesgue-Stieltjes integral.

\noindent {\bf 2010 Mathematics Subject Classification:}\ 49K15, 49J15.

\end{quote}
\section{Introduction}
It is well known that optimal control problems with state constraints are models of importance, but one usually faces with a lot of difficulties in analyzing them. These models have been considered since the early days of the optimal control theory. For instance, the whole Chapter VI of the classical work \cite[pp.~257--316]{Pont_Bolt_Gamk_Mish_1962} is devoted to problems with restricted phase coordinates. There are various forms of the maximum principle for optimal control problems with state constraints; see, e.g.,~\cite{Hartl_Sethi_Vickson_1995}, where the relations between several forms are shown and a series of numerical illustrative examples have been solved.

To deal with state constraints, one has to use functions of bounded variation, Borel measurable functions, Lebesgue-Stieltjes integral, nonnegative measures on the $\sigma-$algebra of the Borel sets, the Riesz Representation Theorem for the space of continuous functions, and so on.

By using the maximum principle presented in \cite[pp.~233--254]{Ioffe_Tihomirov_1979}, Phu \cite{Phu_1989,Phu_1992} has proposed  an ingenious method called \textit{the method of region analysis} to solve several classes of optimal control problems with one state and one control variable, which have both state and control constraints. Minimization problems of the Lagrange type were considered by the author and, among other things, it was assumed that integrand of the objective function is strictly convex with respect to the control variable. To be more precise, the author considered \textit{regular problems}, i.e., the  optimal control problems where the Pontryagin function is strictly convex with respect to the control variable. 

In the present paper, the maximum principle for finite horizon state constrained problems from the book by Vinter \cite[Theorem~9.3.1]{Vinter_2000} is analyzed via parametric examples. The latter has origin in a recent paper by Basco, Cannarsa, and Frankowska \cite[Example~1]{Basco_Cannarsa_Frankowska_2018}, and resembles the optimal growth problem in mathematical economics (see, e.g., \cite[pp.~617--625]{Takayama_1974}). The solution existence of these parametric examples, which are \textit{irregular optimal control problems} in the sense of Phu \cite{Phu_1989,Phu_1992}, is established by invoking Filippov's existence theorem for Mayer problems \cite[Theorem~9.2.i and Section~9.4]{Cesari_1983}. Since the maximum principle is only a necessary condition for local optimal processes, a large amount of additional investigations is needed to obtain a comprehensive synthesis of finitely many processes suspected for being local minimizers. Our analysis not only helps to understand the principle in depth, but also serves as a sample of applying it to meaningful prototypes of  economic optimal growth models. 

Note that the \textit{maximum principle} for finite horizon state constrained problems in \cite[Chapter~9]{Vinter_2000} covers many known ones for smooth problems and allows us to deal with nonsmooth problems by using the \textit{Mordukhovich normal cone} and the \textit{Mordukhovich subdifferential} \cite{Mordukhovich_2006a,Mordukhovich_2006b,Mordukhovich_2018}, which are also called the limiting normal cone and the limiting subdifferential. This principle is a necessary optimality condition which asserts the existence of a nontrivial multipliers set consisting of an \textit{absolutely continuous function}, a \textit{function of bounded variation}, a \textit{Borel measurable function}, and a \textit{real number}, such that the four conditions (i)--(iv) in Theorem~\ref{V_thm9.3.1 necessary condition} below are satisfied. The relationships between these conditions are worthy a detailed analysis. We will present such an analysis via three parametric examples of optimal control problems of the Langrange type, which have five parameters: the first one appears in the description of the objective function, the second one appears in the differential equation, the third one is the initial value, the fourth one is the initial time, and the fifth one is the terminal time. Observe that, in Example~1 of \cite{Basco_Cannarsa_Frankowska_2018}, the terminal time is infinity, the initial value and the initial time are fixed. Problems with unilateral state constraints have been studied in Part 1 (see \cite{VTHuong_Yao_Yen_Part1}) of the paper. Problems with bilateral state constraints are addressed in this Part 2, which is organized as follows. 

Section \ref{Background Materials} presents some background materials including the above-mentioned maximum principle and Filippov's existence theorem for Mayer problems. Control problems with bilateral state constraints are studied in Section \ref{Example 3}. Some concluding remarks are given in Section~\ref{Conclusions}.

In comparison with Part 1, to deal with bilateral state constraints, herein we have to prove a series delicate lemmas and auxiliary propositions. Moreover, the synthesis of finitely many processes suspected for being local minimizers is rather sophisticated, and it requires a lot of refined arguments. 

\section{Background Materials}\label{Background Materials}
In this section, we give some notations, definitions, and results that will be used repeatedly in the sequel. 

\subsection{Notations and Definitions }
 The symbol $\R$ (resp., $\N)$ denotes the set of real numbers (resp., the set of positive integers). The norm in the $n$-dimensional Euclidean space $\R^n$ is denoted by $\|.\|$. For a subset $C\subset\R^n$, we abbreviate its \textit{convex hull} to $\co C$. For a  set-valued map $F: \R^n \rightrightarrows \R^m$,  we call the set $${\rm gph}\, F := \{(x,y) \in \R^n \times \R^m\,:\, y\in F(x) \}$$ the \textit{graph} of $F$.
 
 Let $\Omega \subset \R^n$ be a closed set and $\bar{v} \in \Omega$. The \textit{Fr\'echet} (or \textit{regular}) \textit{normal cone} to
 $\Omega\subset\R^n$ at $\bar{v}$ is given by
 \begin{equation*}
  \widehat
 	N_{\Omega}(\bar{v})=\left\{v'\in\R^n\,:\,
 	\displaystyle\limsup_{v\xrightarrow{\Omega}\bar{v}}\,\displaystyle\frac{\langle
 		v',v-\bar{v} \rangle}{\|v-\bar{v} \|}\leq 0\right\},\end{equation*} where
 $v\xrightarrow{\Omega}\bar{v}$ means $v\to \bar{v}$ with $v\in\Omega$.  The {\it Mordukhovich} (or {\it limiting}) {\it normal cone} to $\Omega$ at $\bar{v}$ is defined by
 \begin{equation*} N_{\Omega}(\bar{v})
 		=\big\{v'\in\R^n\,:\,\exists \mbox{ sequences } v_k\to \bar{v},\ v_k'\rightarrow v' \mbox{ with } v_k'\in \widehat
 		N_{\Omega}(v_k)\; \mbox{for all}\; k\in \N\big\}.\end{equation*}
Given an extended real-valued function $\varphi: \R^n \rightarrow \R\cup\{-\infty, +\infty\}$, one defines the \textit{epigraph} of $\varphi$ by  $\epi \varphi=\{(x, \mu)\in \R^n\times \R \,:\, \mu \geq \varphi (x)\}$. The \textit{Mordukhovich} (or \textit{limiting}) \textit{subdifferential} of $\varphi$ at $\bar x\in \R^n$ with $|\varphi (\bar x)|< \infty$ is defined by
\begin{equation*}
\partial\varphi(\bar x)=\big\{x^*\in \R^n \;:\; (x^*, -1)\in N\big((\bar x, \varphi(\bar x)); \epi \varphi \big)\big\}.
\end{equation*}
If $|\varphi (x)|= \infty$, then one puts $\partial\varphi(\bar x)=\emptyset$.
The reader is referred to \cite[Chapter~1]{Mordukhovich_2006a} and \cite[Chapter~1]{Mordukhovich_2018} for comprehensive treatments of the Fr\'echet normal cone, the limiting normal cone, the limiting subdifferential, and the related calculus rules.
 
For a given segment $[t_0, T]$ of the real line, we denote the $\sigma$-algebra of its Lebesgue measurable subsets (resp., the $\sigma$-algebra of its Borel sets) by $\mathcal{L}$ (resp.,~$\mathcal{B}$). The Sobolev space  $W^{1,1}([t_0, T], \R^n)$ is the linear space of the absolutely continuous functions $x:[t_0, T] \to \R^n$ endowed with the norm $$\|x\|_{W^{1,1}}=\|x(t_0)\|+\displaystyle\int_{t_0}^T \|\dot x(t)\| dt$$ (see, e.g., \cite[p.~21]{Kolmogorov_Fomin_1970} for this and another equivalent norm).

As in \cite[p.~321]{Vinter_2000}, we consider the following \textit{finite horizon optimal control problem of the Mayer type}, denoted by $\mathcal M$,
\begin{equation}\label{cost functional_FP}
\mbox{Minimize}\ \; g(x(t_0), x(T)),
\end{equation}
over $x \in W^{1,1}([t_0, T], \R^n)$  and measurable functions $u:[t_0, T] \to \R^m$ satisfying
\begin{equation}\label{state control system_FP}
	\begin{cases}
		\dot x(t)=f(t, x(t), u(t)),\quad &\mbox{a.e.\ } t\in [t_0, T]\\
		(x(t_0), x(T))\in C\\
		u(t)\in U(t), &\mbox{a.e.\ } t\in [t_0, T]\\
		h(t, x(t))\leq 0, & \forall t\in [t_0, T],
	\end{cases}
\end{equation}
where $[t_0, T]$ is a given interval, $g: \R^n\times \R^n \to \R$, $f: [t_0, T]\times \R^n\times \R^m \to \R^n$, and $h:[t_0, T]\times \R^n \to \R$ are given functions, $C\subset \R^n\times \R^n$ is a closed set, and $U: [t_0, T]\rightrightarrows \R^m$ is a set-valued map. 

A measurable function $u:[t_0, T] \to \R^m$ satisfying $u(t)\in U(t)$ a.e. $t\in [t_0, T]$ is called a \textit{control function}. A \textit{process} $(x, u)$ consists of a control function $u$ and an arc $x \in W^{1,1}([t_0, T]; \R^n)$ that is a solution to the differential equation in \eqref{state control system_FP}. A \textit{state trajectory} $x$ is the first component of some process $(x, u)$. A process $(x, u)$ is called \textit{feasible} if the state trajectory satisfies the \textit{endpoint constraint} $(x(t_0), x(T))\in C$ and the \textit{state constraint} $h(t, x(t))\leq 0$ for all $t\in [t_0, T]$.

Due to the appearance of the state constraint, the problem $\mathcal M$ in \eqref{cost functional_FP}--\eqref{state control system_FP} is said to be an \textit{optimal control problem with state constraints}. But, if the inequality $h(t, x(t))\leq 0$ is fulfilled for every $(t, x(t))$ with $t \in [t_0, T]$ and $x \in W^{1,1}([t_0, T]; \R^n)$ (for example, when $h$ is constant function having a fixed nonpositive value), i.e., the condition $h(t, x(t))\leq 0$ for all  $t\in [t_0, T]$ can be removed from~\eqref{state control system_FP}, then one says that $\mathcal M$ an \textit{optimal control problem without state constraints}.

 The \textit{Hamiltonian} $\mathcal H: [t_0, T]\times \R^n \times \R^n\times \R^m \to \R$ of \eqref{state control system_FP} is defined by 
\begin{equation}\label{Hamiltonian}
\mathcal H(t, x, p, u):=p.f(t, x, u)=\displaystyle\sum_{i=1}^np_if_i(t, x, u). 
\end{equation}

\begin{definition}\label{local_minimizer} {\rm 	
	A feasible process $(\bar x, \bar u)$ is called a $W^{1,1}$ \textit{local minimizer} for  $\mathcal M$ if there exists $\delta>0$ such that $g(\bar x(t_0), \bar x(T))\leq g(x(t_0), x(T))$ for any feasible processes $(x, u)$ satisfying $\|\bar x-x\|_{W^{1,1}} \leq \delta$.}
	\end{definition}
	
	\begin{definition}\label{global_minimizer} {\rm 	
			A feasible process $(\bar x, \bar u)$ is called a $W^{1,1}$ \textit{global minimizer} for  $\mathcal M$ if, for any feasible processes $(x, u)$, one has $g(\bar x(t_0), \bar x(T))\leq g(x(t_0), x(T))$.}
	\end{definition}
	
\begin{definition}[{See \cite[p.~329]{Vinter_2000}}]\label{partial_hybrid subdiff} {\rm 
The \textit{partial hybrid subdifferential} $\partial^>_x h(t, x)$ of $h(t,x)$ w.r.t. $x$ is given by
\begin{align}\label{h-partial subdiff}
\partial^>_x h(t, x):=\co\big\{\xi \,:\, &\mbox{ there exists }(t_i, x_i)\overset{h}{\rightarrow}(t, x) \mbox{ such that }\nonumber\\
&\ h(t_k, x_k)>0 \mbox{ for all } k \mbox{ and }\nabla_xh(t_k, x_k)\to \xi\big\}, 
\end{align} where $(t_k, x_k)\overset{h}{\rightarrow}(t, x)$ means that $(t_k, x_k)\rightarrow (t, x)$ and $h(t_k, x_k)\rightarrow h(t, x)$ as $k\to\infty$.}
\end{definition}

\subsection{A Maximum Principle for State Constrained Problems}

Due to the appearance of the state constraint $h(t, x(t))\leq 0$ in $\mathcal{M}$, one has to introduce a multiplier that is an element in the topological dual $C^*([t_0, T]; \R)$ of the space of continuous functions $C([t_0, T]; \R)$ with the supremum norm. By the Riesz Representation Theorem (see, e.g., \cite[Theorem~6, p.~374]{Kolmogorov_Fomin_1970} and \cite[Theorem~1, pp.~113--115]{Luenberger_1969}), any bounded linear 
functional $f$ on $C([t_0, T]; \R)$ can be uniquely represented in the form 
\begin{equation*}
f(x) =\int_{[t_0,T]} x(t) dv(t),
\end{equation*}
where $v$ is \textit{a function of bounded variation} on $[t_0, T]$ which vanishes at $t_0$ and which are continuous from the right at every point $\tau\in (t_0, T)$, and $\displaystyle\int_{[t_0,T]} x(t) dv(t)$ is the Riemann-Stieltjes integral of $x$ with respect to $v$ (see, e.g., \cite[p.~364]{Kolmogorov_Fomin_1970}). 
The set of the elements of $C^*([t_0, T]; \R)$ which are given by nondecreasing functions $v$ is denoted by $C^\oplus(t_0, T)$.

Every $v \in C^*([t_0, T]; \R)$ corresponds to \textit{a finite regular measure}, denoted by $\mu_v$, on the $\sigma$-algebra ${\mathcal B}$ of the Borel subsets of $[t_0, T]$ by the formula
	\begin{equation*}
	\mu_v(A) :=\int_{[t_0,T]} \chi_A(t) dv(t),
	\end{equation*} where $\chi_A(t)=1$ for $t\in A$ and $\chi_A(t)=0$ if $t\notin A$. Due to the correspondence $v\mapsto\mu_v$, we call every element $v\in C^*([t_0, T]; \R)$ a ``measure" and identify $v$ with $\mu_v$. Clearly,  the measure corresponding to each $v\in C^\oplus(t_0, T)$ is nonnegative.
	
	The integrals $\displaystyle\int_{[t_0, t)}\nu(s)d\mu(s)$ and $\displaystyle\int_{[t_0, T]}\nu(s)d\mu(s)$ of a Borel measurable function $\nu$ in next theorem are understood in the sense of the Lebesgue-Stieltjes integration \cite[p.~364]{Kolmogorov_Fomin_1970}. 
	
The $\sigma$-algebra of the Borel sets in $\R^m$ is denoted by $\mathcal{B}^m$.

\begin{theorem}[{See \cite[Theorem~9.3.1]{Vinter_2000}}]\label{V_thm9.3.1 necessary condition}
Let $(\bar x, \bar u)$ be a $W^{1,1}$ local minimizer for $\mathcal M$. Assume that for some $\delta>0$, the following hypotheses are satisfied:
\begin{enumerate}[\rm (H1)]
\item $f(., x, .)$ is $\mathcal{L}\times \mathcal{B}^m$ measurable, for fixed $x$. There exists a Borel measurable function $k(., .):[t_0, T]\times \R^m \to \R$ such that $t \mapsto k(t, \bar u(t))$ is integrable and 
\begin{equation*}
\| f(t, x, u)-f(t, x', u)\|\leq k(t, u)\|x-x'\|, \quad \forall x, x'\in \bar x(t)+\delta\bar B,\; \forall u \in U(t), \mbox{a.e.;}
\end{equation*}
\item $\gph U$ is a Borel set in $[t_0,T]\times\R^m$;
\item $g$ is Lipschitz continuous on the ball $(\bar x(t_0), \bar x(T))+\delta\bar B$;
\item $h$ is upper semicontinuous and there exists $K>0$ such that 
\begin{equation*}
\| h(t, x)-h(t, x')\|\leq K\|x-x'\|, \quad \forall x, x'\in \bar x(t)+\delta\bar B,\; \forall t \in [t_0, T].
\end{equation*}
\end{enumerate}
Then there exist $p\in W^{1,1}([t_0, T]; \R^n)$, $\gamma \geq 0$, $\mu \in C^\oplus(t_0, T)$, and a Borel measurable  function $\nu:[t_0, T]\to  \R^n$ such that $(p, \mu, \gamma)\neq (0, 0, 0)$, and for $q(t):=p(t)+\eta(t)$ with $\eta(t):=
\displaystyle\int_{[t_0, t)}\nu(s)d\mu(s)$ if $t\in [t_0, T)$ and $\eta(T):=\displaystyle\int_{[t_0, T]}\nu(s)d\mu(s)$, the following holds true:
\begin{enumerate}[\rm (i)]
\item $\nu(t)\in \partial^>_x h(t, \bar x(t))\ \mu-\mbox{a.e.};$
 \item $-\dot p(t)\in \co \partial_x\mathcal{H}(t, \bar x(t), q(t), \bar u(t))$ a.e.;
\item $(p(t_0), -q(T))\in \gamma \partial g(\bar x(t_0), \bar x(T))+N_C(\bar x(t_0), \bar x(T))$;
\item $\mathcal{H}(t, \bar x(t), q(t), \bar u(t))=\max_{u\in U(t)}\mathcal{H}(t, \bar x(t), q(t), u)$ a.e.
\end{enumerate}
\end{theorem}

\subsection{Solution Existence in State Constrained Optimal Control}
To recall a solution existence theorem for optimal control problems with state constraints of the Mayer type, we will use the notations and concepts given in \cite[Section~9.2]{Cesari_1983}. Let $A$ be a subset of $\R\times \R^n$ and $U: A\rightrightarrows \R^m$ be a set-valued map defined on $A$. Let $$M:=\{(t, x, u)\in \R\times\R^n\times \R^m \;:\; (t, x)\in A,\ u \in U(t, x)\},$$ and $f=(f_1, f_2, \dots, f_n): M \to \R^n$ be a single-valued map defined on $M$. Let $B$ be a given subset of $\R\times \R^n\times\R\times \R^n$ and $g: B \to \R$ be a real function defined on $B$. Consider the optimal control problem of the Mayer type
\begin{equation}\label{cost functional_SET}
\mbox{Minimize}\ \; g(t_0, x(t_0), T,  x(T))
\end{equation} over $x \in W^{1,1}([t_0, T]; \R^n)$  and measurable functions $u:[t_0, T]~\to~\R^m$ satisfying
\begin{equation}\label{state control system_SET}
	\begin{cases}
		\dot x(t)=f(t, x(t), u(t)),\quad &\mbox{a.e.\ } t\in [t_0, T]\\
		(t, x(t))\in A, & \mbox{for all }t\in [t_0, T]\\
		(t_0, x(t_0), T, x(T))\in B\\
		u(t)\in U(t, x(t)), &\mbox{a.e.\ } t\in [t_0, T],
	\end{cases}
\end{equation}
where $[t_0, T]$ is a given interval. The problem \eqref{cost functional_SET}--\eqref{state control system_SET} will be denoted by $\mathcal{M}_1$.

A \textit{feasible process} for $\mathcal M_1$ is a pair of functions $(x, u)$ with $x: [t_0, T] \to \R^n$ being absolutely continuous on $[t_0, T]$, $u:[t_0, T] \to \R^m$ being measurable, such that  all the requirements in \eqref{state control system_SET} are satisfied. If  $(x, u)$ is a feasible process for $\mathcal M_1$, then $x$
is said to be a \textit{feasible trajectory}, and $u$ a \textit{feasible control function} for $\mathcal M_1$. The set of all feasible processes for $\mathcal M_1$ is denoted by $\Omega$.

Let $A_0=\big\{t\in\mathbb R\,:\, \exists x\in \mathbb R^n\ {\rm s.t.}\ (t,x)\in A\big\}$, i.e., $A_0$ is the projection of $A$ on the $t-$axis. Set $$A(t)=\big\{x\in \R^n \;:\; (t, x)\in A\big\}\quad\; (t\in A_0)$$ and  $$Q(t, x)=\big\{z\in \R^n \;:\; z=f(t, x, u),\ u \in U(t, x)\big\} \quad\; ((t, x)\in A).$$

The forthcoming statement is called \textit{Filippov's Existence Theorem for Mayer problems}.
\begin{theorem}[{see \cite[Theorem~9.2.i and Section~9.4]{Cesari_1983}}]\label{Filippov's_Existence_Thm}
Suppose that $\Omega$ is nonempty, $B$ is closed, $g$ is lower semicontinuous on $B$, $f$ is continuous on $M$ and, for almost every $t\in [t_0, T]$, the sets $Q(t, x)$, $x\in A(t)$, are convex. Moreover,  assume either that $A$ and $M$ are compact or that $A$ is not compact but closed and the following three conditions hold
\begin{enumerate}[\rm (a)]
\item For any $\varepsilon\geq 0$, the set $M_\varepsilon:=\{(t, x, u)\in M \;:\; \|x\| \leq \varepsilon\}$ is compact;
\item There is a compact subset $P$ of $A$ such that every feasible trajectory $x$ of $\mathcal M_1$ passes through at least one point of $P$;
\item There exists $c\geq 0$ such that 
\begin{equation*}
x_1 f_1(t, x, u) + x_2 f_2(t, x, u)+\dots+x_n f_n(t, x, u) \leq c (\|x\|^2+1)\quad\; \forall (t, x, u)\in M.
\end{equation*}
\end{enumerate}
Then, $\mathcal M_1$ has a $W^{1, 1}$ global minimizer.
\end{theorem}
Clearly, condition (b) is satisfied if the initial point $(t_0, x(t_0))$ or the end point $(T, x(T))$ is fixed. As shown in \cite[p.~317]{Cesari_1983}, the following condition implies (c): 
\begin{enumerate}[\rm ($c_0$)]
\item \textit{There exists $c\geq 0$ such that $
\|f(t, x, u)\| \leq c (\|x\|+1)$ for all $(t, x, u)\in M$.}
\end{enumerate}

	    	\section{Optimal Control Problems with Bilateral State Constraints}\label{Example 3}
	    	 By $(FP_3)$ we denote the finite horizon optimal control problem of the Lagrange type
	    		\begin{equation} \label{cost functional_FP_3}
	    		\mbox{Minimize}\ \; J(x,u)=\int_{t_0}^{T} \big[-e^{-\lambda t}(x(t)+u(t))\big] dt
	    		\end{equation}
	    		over $x \in W^{1,1}([t_0, T], \R)$  and measurable functions $u:[t_0, T] \to\R$ satisfying
	    		\begin{equation} \label{state control system_FP_3}
	    		\begin{cases}
	    		\dot x(t)=-au(t),\quad &\mbox{a.e.\ } t\in [t_0, T]\\
	    		x(t_0)=x_0\\
	    		u(t)\in [-1, 1], &\mbox{a.e.\ } t\in [t_0, T]\\
	    	   -1\leq x(t)\leq 1, & \forall t\in [t_0, T]
	    		\end{cases}
	    		\end{equation}
	    		with $a > \lambda>0$, $T>t_0\geq 0$, and  $-1\leq x_0 \leq 1$ being given.
	    		
	    		\medskip
	    	To treat  $(FP_3)$  in \eqref{cost functional_FP_3}--\eqref{state control system_FP_3} as a problem of the Mayer type, we set $x(t)=(x_1(t), x_2(t))$, where $x_1(t)$ plays the role of $x(t)$ in $(FP_3)$ and \begin{equation}\label{phase_second_component_FP3a}x_2(t):= \int_{t_0}^{t} \big[-e^{-\lambda\tau}(x_1(\tau)+u(\tau))\big] d\tau\end{equation} for all $t\in [0, T]$. Thus, $(FP_3)$ is equivalent to the problem
	    		\begin{equation}\label{cost functional_FP_3a}
	    		\mbox{Minimize}\ \; x_2(T)
	    		\end{equation}
	    		over $x=(x_1, x_2) \in W^{1,1}([t_0, T], \R^2)$  and measurable functions $u:[t_0, T] \to\R$ satisfying
	    		\begin{equation} \label{state control system_FP_3a}
	    		\begin{cases}
	    		\dot x_1(t)=-au(t),\quad &\mbox{a.e.\ } t\in [t_0, T]\\
	    		\dot x_2(t)=-e^{-\lambda t}(x_1(t)+u(t)), &\mbox{a.e.\ } t\in [t_0, T]\\
	    		(x(t_0), x(T))\in \{(x_0, 0)\}\times \R^2\\
	    		u(t)\in [-1, 1], &\mbox{a.e.\ } t\in [t_0, T]\\
	    		-1\leq x_1(t)\leq 1,\ &\forall t\in [t_0, T].
	    		\end{cases}
	    		\end{equation} The problem \eqref{cost functional_FP_3a}--\eqref{state control system_FP_3a} is abbreviated to  $(FP_{3a})$.
	    				
	    	\subsection{Solution Existence}\label{SE_2-sided_S_constraint}
	    		To verify that $(FP_{3a})$ is of the form $\mathcal{M}_1$ (see Subsection 2.3), one can choose $n=2$, $m=1$, $A=[t_0, T]~\times[-1, 1]\times \R$, $U(t, x)=[-1, 1]$ for all $(t, x)\in A$, $B=\{t_0\}\times\{(x_0, 0)\}\times \R\times \R^2$, $g(t_0, x(t_0), T, x(T))=x_2(T)$, $M=A\times [-1, 1]$, $f(t, x, u)=(-au, -e^{-\lambda t}(x_1+u))$ for all $(t, x, u)\in M$. To show that $(FP_{3a})$ satisfies all the assumptions of Theorem~\ref{Filippov's_Existence_Thm}, we can use the arguments given in Subsection~3.1 in Part 1~(\cite{VTHuong_Yao_Yen_Part1}), except those related to the convexity of the sets $Q(t,x)$ and the compactness of $M_\varepsilon$.
	    				
	    	By the formula for $A$, one has  $A_0=[t_0, T]$ and $A(t)=[-1, 1]\times \R$ for all $t\in A_0$. Thus, the requirement in Theorem~\ref{Filippov's_Existence_Thm} on the convexity of the sets  $Q(t, x)$, $x\in A(t)$, for almost every $t\in [t_0, T]$ is satisfied. Since $M=[t_0, T]~\times[-1, 1]\times \R\times [-1, 1]$, for any $\varepsilon\geq 0$, one has the expression
	    	$M_\varepsilon=\{(t, x, u)\in[t_0, T]~\times[-1, 1]\times \R\times [-1, 1] \,:\, \|x\| \leq \varepsilon\},$ which justifies the compactness of	$M_\varepsilon$.

	    		 Theorem~\ref{Filippov's_Existence_Thm} tells us that $(FP_{3a})$ has a $W^{1, 1}$ global minimizer. Thus, by the equivalence of~$(FP_{3})$ and $(FP_{3a})$, we can assert that $(FP_{3})$ has a $W^{1, 1}$ global minimizer.
	    	    				
	    		\subsection{Necessary Optimality Conditions} 
	    		 To solve problem $(FP_3)$ by applying Theorem~\ref{V_thm9.3.1 necessary condition}, note that $(FP_{3a})$ is in the form of $\mathcal M$ with  $g(x, y)=y_2$, $$f(t, x, u)=(-au, -e^{-\lambda t}(x_1+u)),$$ $C=\{(x_0, 0)\}\times \R^2$, $U(t)=[-1, 1]$, and $h(t, x)=|x_1|-1$ for all $t\in [t_0, T]$, $x=(x_1, x_2) \in \R^2$, $y=(y_1, y_2) \in \R^2$ and $u\in \R$. According to \eqref{Hamiltonian}, the Hamiltonian of $(FP_{3a})$ is the function
	    	    \begin{equation}\label{Hamiltonian_FP_3a}
	    	    \mathcal{H}(t, x, p, u)=-aup_1-e^{-\lambda t}(x_1+u)p_2 \quad \forall(t, x, p, u)\in [t_0, T]\times \R^2\times \R^2 \times \R.
	    	    \end{equation}
	    	    By \eqref{h-partial subdiff}, the partial hybrid subdifferential of $h$ at $(t, x)\in [t_0, T]\times \R^2$ is the set
	    		\begin{equation}\label{hybrid subdiff_FP_3a}
	    		\partial^>_x h(t, x)=
	    		\begin{cases}
	    		\{(-1,0)\}, & \quad \mbox{if}\ x_1\leq -1\\
	    		\emptyset, &\quad \mbox{if} \ |x_1|<1\\
	    		\{(1,0)\}, & \quad \mbox{if}\ x_1\geq 1.\\
	    	    \end{cases}
	    		\end{equation}
	    				
	    		Let $(\bar x, \bar u)$ be a $W^{1,1}$ local minimizer for $(FP_{3a})$. Since the  assumptions (H1)--(H4) of Theorem~\ref{V_thm9.3.1 necessary condition} are satisfied for $(FP_{3a})$, by that theorem one can find $p\in W^{1,1}([t_0, T]; \R^2)$, $\gamma \geq 0$, $\mu \in C^\oplus(t_0, T)$, and a Borel measurable  function $\nu:[t_0, T]\to  \R^2$ such that $(p, \mu, \gamma)\neq (0, 0, 0)$, and for $q(t):=p(t)+\eta(t)$ with \begin{equation*}\label{1st_formula_for_eta_FP3a}\eta(t):=
	    	    \displaystyle\int_{[t_0, t)}\nu(\tau)d\mu(\tau)\quad\ (\forall t\in [t_0, T))\end{equation*} and
	    		\begin{equation*}\label{2nd_formula_for_eta_FP3a}
	    		\eta(T):=\displaystyle\int_{[t_0, T]}\nu(\tau)d\mu(\tau),
	    		\end{equation*}
	    		conditions (i)--(iv) in Theorem~\ref{V_thm9.3.1 necessary condition} hold true.
	    		
	            \textbf{Condition (i)}: Note that
	    		\begin{align*}
	    		&  \mu \{t\in [t_0, T] \,:\, \nu(t) \notin \partial^>_x h(t, \bar x(t))\}\\ &=\ \mu \{t\in [t_0, T] \,:\, \partial^>_x h(t, \bar x(t))=\emptyset\}+\mu \{t\in [t_0, T] \,:\, \partial^>_x h(t, \bar x(t))\neq\emptyset,\; \nu(t) \notin \partial^>_x h(t, \bar x(t))\}.
	    		\end{align*}
	    		Since $-1\leq x_1(t)\leq 1$ for every $t$, combining this with \eqref{hybrid subdiff_FP_3a} gives
	    		\begin{align*}
	    		\mu \{t\in [t_0, T] \,:\, \nu(t) \notin \partial^>_x h(t, \bar x(t))\} =\ &\mu \{t\in [t_0, T] \,:\, -1< \bar x_1(t)< 1\}\\ 
	    		&+\mu \{t\in [t_0, T] \,:\, \bar x_1(t)=1,\; \nu(t) \neq (1, 0)\}\\
	    		&+ \mu \{t\in [t_0, T] \,:\, \bar x_1(t)=-1,\; \nu(t) \neq (-1, 0)\}.	    	
	    		\end{align*}
	    	    So, from (i)  it follows that 
	    		\begin{equation}\label{1st_condition_for_mu_FP3a}\mu \{t\in [t_0, T] \,:\,-1< \bar x_1(t)< 1\}=0,\end{equation}
	    			\begin{equation}\label{2nd_condition_for_mu_FP3a}\mu \big\{t\in [t_0, T] \,:\, \bar x_1(t)=1,\; \nu(t) \neq (1, 0)\big\}=0,\end{equation}
	    			 \begin{equation}\label{3th_condition_for_mu_FP3a}\mu \big\{t\in [t_0, T] \,:\, \bar x_1(t)=-1,\; \nu(t) \neq (-1, 0)\big\}=0.\end{equation}
	    	
	    		\textbf{Condition (ii)}: By \eqref{Hamiltonian_FP_3a}, $\mathcal{H}$ is differentiable in $x$ and $\partial_x\mathcal{H}(t, x, p, u)=\{(-e^{-\lambda t}p_2, 0)\}$ for all $(t, x, p, u)\in [t_0, T]\times \R^2\times \R^2 \times \R$. Thus, (ii) implies that $-\dot p(t) =(-e^{-\lambda t}q_2(t), 0)$ for a.e. $t\in [t_0, T]$. Hence, $\dot p_1(t) = e^{-\lambda t}q_2(t)$  for a.e. $t\in [t_0, T]$ and $p_2(t)$ is  a constant for all $t\in [t_0, T]$.
	    				
	    		\textbf{Condition (iii)}: By the formulas for $g$ and $C$, $\partial g(\bar x(t_0), \bar x(T))=\{(0, 0, 0, 1)\}$ and $N_C(\bar x(t_0), \bar x(T))=\R^2\times\{(0, 0)\}$. Thus, (iii) yields
	    		$$(p(t_0), -q(T))\in \{(0, 0, 0, \gamma)\}+\R^2\times\{(0, 0)\},$$ which means that $q_1(T)=0$ and $q_2(T)=-\gamma$.	
	    		
	    	    \textbf{Condition (iv)}: By \eqref{Hamiltonian_FP_3a}, from (iv) one gets
	    		\begin{equation*}
	    		-a\bar u(t)q_1(t)-e^{-\lambda t}[\bar x_1(t)+\bar u(t)]q_2(t)=\max_{u\in [-1, 1]}\left\{-auq_1(t)-e^{-\lambda t}[\bar x_1(t)+u]q_2(t) \right\}\ \mbox{a.e.}\; t\in [t_0,T]
	    		\end{equation*}
	    		or, equivalently,
	    		\begin{equation*}\label{min_condition_FP_3a}
	    		[aq_1(t)+e^{-\lambda t}q_2(t)]\bar u(t)=\min_{u\in [-1, 1]}\left\{[aq_1(t)+e^{-\lambda t}q_2(t)]u\right\}\ \mbox{a.e.}\; t\in [t_0,T].
	    		\end{equation*}	
	    		
	    		If the curve $\bar x_1(t)$ remains in the interior of the domain $[-1,1]$
	 for all $t$ from an open interval $(\tau_1,\tau_2)$ of the time axis and touches the boundary of the domain at the moments $\tau_1$ and $\tau_2$, then it must have some special form. A formal formulation of this observation is as follows.
	    			    		    			    	
	    			    	\begin{Proposition}\label{Lemma_main_tool}
	    			    		Suppose that $[\tau_1, \tau_2]$, $\tau_1<\tau_2$, is a subsegment of $[t_0, T]$ with $\bar x_1(t)\in (-1, 1)$ for all $t\in (\tau_1, \tau_2)$. Then, next statements hold true.
	    			    		\begin{enumerate}[\rm S1)]
	    			    			\item If $\bar x_1(\tau_1)=-1$ and $\bar x_1(\tau_2)=1$, then $\tau_2-\tau_1 = 2a^{-1}$ and
	    			    			\begin{equation*}
	    			    			\bar x_1(t)=-1+a(t-\tau_1), \quad  t\in [\tau_1, \tau_2].
	    			    			\end{equation*}
	    			    			\item If $\bar x_1(\tau_1)=1$ and $\bar x_1(\tau_2)=-1$, then $\tau_2-\tau_1 = 2a^{-1}$ and
	    			    			\begin{equation*}
	    			    			\bar x_1(t)=1-a(t-\tau_1), \quad  t\in [\tau_1, \tau_2].
	    			    			\end{equation*}
	    			    			\item If $\bar x_1(\tau_1)=\bar x_1(\tau_2)=-1$, then $\tau_2-\tau_1< 4a^{-1}$ and
	    			    			\begin{equation*}
	    			    			\bar x_1(t)=
	    			    			\begin{cases}
	    			    			-1+a(t-\tau_1), \quad & t\in [\tau_1, \hat t]\\
	    			    			-1-a(t-\tau_2), &  t\in (\hat t, \tau_2],
	    			    			\end{cases}
	    			    			\end{equation*}
	    			    			where  $\hat t:=(\tau_1+\tau_2)/2$.
	    			    			\item The situation where $\bar x_1(\tau_1)=\bar x_1(\tau_2)=1$ cannot happen.
	    			    		\end{enumerate}
	    			    	\end{Proposition}
	    			    	\begin{proof}
	    			    		Choose $\varepsilon_1>$ and $\varepsilon_2>0$ small enough so as $[\tau_1+\varepsilon_1, \tau_2-\varepsilon_2]\subset [\tau_1, \tau_2]$. Then, $\bar x_1(t)\in (-1, 1)$ for all $ t\in [\tau_1+\varepsilon_1, \tau_2-\varepsilon_2]$, i.e., $h(t, \bar x(t))<0$ for all $t \in [\tau_1+\varepsilon_1, \tau_2-\varepsilon_2]$. Thus, applying Proposition~4.3 in Part 1~(\cite{VTHuong_Yao_Yen_Part1}) with $(FP_{3a})$ in the place of $(FP_{2a})$ in its formulation, one finds that the formula for $\bar x_1(.)$ on $[\tau_1+\varepsilon_1, \tau_2-\varepsilon_2]$ belongs to one of the following categories C1$-$C3:
	    			    		\begin{equation*}\label{interior_traj_1_FP3a}
	    			    		\bar x_1(t)=\bar x_1(\tau_1+\varepsilon_1)+a(t-\tau_1-\varepsilon_1), \quad t\in [\tau_1+\varepsilon_1, \tau_2-\varepsilon_2],
	    			    		\end{equation*}
	    			    		\begin{equation*}\label{interior_traj_2_FP3a}
	    			    		\bar x_1(t)=\bar x_1(\tau_1+\varepsilon_1)-a(t-\tau_1-\varepsilon_1), \quad t\in [\tau_1+\varepsilon_1, \tau_2-\varepsilon_2],
	    			    		\end{equation*}
	    			    		and
	    			    		\begin{equation*} \label{interior_traj_FP3a}
	    			    		\bar x_1(t)=
	    			    		\begin{cases}
	    			    		\bar x_1(\tau_1+\varepsilon_1)+a(t-\tau_1-\varepsilon_1), \quad & t\in [\tau_1+\varepsilon_1, t_{\zeta}]\\
	    			    		\bar x_1(t_{\zeta})-a(t-t_\zeta), & t\in (t_{\zeta}, \tau_2-\varepsilon_2],
	    			    		\end{cases}
	    			    		\end{equation*}
	    			    		where $t_{\zeta}$ is some point in $(\tau_1+\varepsilon_1, \tau_2-\varepsilon_2)$. 
	    			    		
	    			    		To prove the statement S1, let $\varepsilon_2=k^{-1}$ with $k$ being a positive integer, as large as  $k^{-1}\in (\tau_1+\varepsilon_1, \tau_2)$. Since for each $k$ the formula for $\bar x_1(.)$ on $[\tau_1+\varepsilon_1,\tau_2-k^{-1}]$ must be of  the three types C1--C3, by the Dirichlet principle there must exist a subsequence $\{k'\}$ of $\{k\}$ such that the corresponding formulas belong to a fixed category. If the latter is happens to be C2, then by the continuity of $\bar x_1(.)$ one has $$\bar x_1(\tau_2)=\displaystyle\lim_{k'\to\infty} \bar x_1(\tau_2-\dfrac{1}{k'})=\displaystyle\lim_{k'\to\infty} \Big[\bar x_1(\tau_1+\varepsilon_1)-a(\tau_2-\dfrac{1}{k'}-\tau_1-\varepsilon_1)\Big]=\bar x_1(\tau_1+\varepsilon_1)-a(\tau_2-\tau_1-\varepsilon_1).$$ This is impossible, because $\bar x_1(\tau_2)=1$. Similarly, the situation where the fixed category is~C3 must also be excluded. In the case where the formulas for $\bar x_1(.)$ belong to the category~C1, we have
	    			    		\begin{equation*}
	    			    		\bar x_1(t)=\bar x_1(\tau_1+\varepsilon_1)+a(t-\tau_1-\varepsilon_1), \quad  t\in [\tau_1+\varepsilon_1, \tau_2].
	    			    		\end{equation*} Now, letting $\varepsilon_1$ tend to zero and using continuity of $\bar x_1(.)$, we obtain 
	    			    		\begin{equation*}
	    			    		\bar x_1(t)=\bar x_1(\tau_1)+a(t-\tau_1), \quad  t\in [\tau_1, \tau_2].
	    			    		\end{equation*}
	    			    		As $\bar x_1(\tau_1)=-1$, the statement S1 is proved.
	    			    		
	    			    		The statements S2 and S3 are proved similarly.
	    			    		
	    			    		To prove the assertion S4, it suffices to apply the arguments of the second part of the analysis of Subcase 4b in Subsection~4.2 in Part 1~(\cite{VTHuong_Yao_Yen_Part1}).
	    			    	\end{proof}	
	    			    		
	    			    		  The forthcoming technical lemma will be in use very frequently.
	    			    		 	    			    		
	    			    		\begin{Lemma}\label{lemma_technical1}
	    		 	Given any $t_1, t_2\in [t_0,T]$, $t_1<t_2$, one puts 
	    		 	\begin{equation}\label{integral_on_segment}
	    		 	J(x, u)|_{[t_1, t_2]}:=\displaystyle\int_{t_1}^{t_2} \big[-e^{-\lambda t}\big(x_1(t)+u(t)\big)\big] dt
	    		 	\end{equation}
	    		 	for any feasible process $(x, u)$ of $(FP_{3a})$.  If $(\widetilde x, \widetilde u)$ and $(\check x, \check u)$ are feasible processes for $(FP_{3a})$ with $\widetilde x_1(t)=1$ for all $t\in [t_1, t_2]$ and 
	    		 	\begin{equation}\label{hat_x1}
	    		 	\check x_1(t)=
	    		 	\begin{cases}
	    		 	1-a(t-t_1), &\quad t\in [t_1, \check t ]\\
	    		 	1+a(t-t_2), &\quad t\in (\check t, t_2],
	    		 	\end{cases}
	    		 	\end{equation}
	    		 	where $\check t:=2^{-1}(t_1+t_2)$, then one has 
	    		 	\begin{eqnarray}\label{identity}
	    		 		J(\check x, \check u)|_{[t_1, t_2]}-J(\widetilde x, \widetilde u)|_{[t_1, t_2]}=\dfrac{1}{\lambda}\big(\dfrac{a}{\lambda}-1\big)\Delta(t_1, t_2)
	    		 	\end{eqnarray} with \begin{eqnarray}\label{Delta}\Delta(t_1, t_2):=e^{-\lambda t_1}-2e^{-\frac{1}{2}\lambda (t_1+t_2)}+e^{-\lambda t_2}.	\end{eqnarray} Besides, it holds that $\Delta(t_1, t_2)>0$ and $J(\check x, \check u)|_{[t_1, t_2]}>J(\widetilde x, \widetilde u)|_{[t_1, t_2]}$.
	    		 \end{Lemma}
	    		 \begin{proof} Using the equation $\dot x_1(t)=-au(t)$ in \eqref{state control system_FP_3a}, which is fulfilled for almost all $t\in [t_0, T]$, and the assumed properties of the processes  $( \widetilde x, \widetilde u)$ and $(\check x, \check u)$, we have $\widetilde u(t)=0$ for almost all $t\in [t_1, t_2]$ and 
	    		 	\begin{equation*}
	    		 	\check u(t)=
	    		 	\begin{cases}
	    		 	1, &\quad \mbox{a.e.}\ \, t\in [t_1, \check t ]\\
	    		 	-1, &\quad \mbox{a.e.}\ \, t\in (\check t, t_2].
	    		 	\end{cases}
	    		 	\end{equation*}
	    		 	Since $\check x (\cdot)$ is a feasible trajectory for $(FP_{3a})$, one has $\check x(\check t)\geq -1$, i.e., $t_2-t_1\leq 4a^{-1}$. 
	    		 	
	    		 By the formulas for $\widetilde x_1$ and $\widetilde u$ on $[t_1, t_2]$, 
	    		 	\begin{equation*}
	    		 	J(\widetilde x, \widetilde u)|_{[t_1, t_2]}=\int_{t_1}^{t_2} \big[-e^{-\lambda t}\big(\widetilde x_1(t)+\widetilde u(t)\big)\big] dt=-\int_{t_1}^{t_2} e^{-\lambda t}dt=\dfrac{1}{\lambda}e^{-\lambda t_2}-\dfrac{1}{\lambda}e^{-\lambda t_1}.
	    		 	\end{equation*}
	    		 	Similarly, from the formulas for $\check x_1$ and $\check u$ on $[t_1, t_2]$ it follows that
	    		 	\begin{align*}
	    		 	J(\check x, \check u)|_{[t_1, t_2]}=&\int_{t_1}^{t_2} \big[-e^{-\lambda t}\big(\check x_1(t)+\check u(t)\big)\big] dt\\
	    		 	=& \int_{t_1}^{\check t}\big[-e^{-\lambda t}\big((1-a(t-t_1))+1 \big)\big] dt\\& +\int_{\check t}^{t_2} \big[-e^{-\lambda t}\big((1+a(t-t_2)) -1  \big)\big] dt\\
	    		 	=& \int_{t_1}^{\check t}e^{-\lambda t}\big[a(t-t_1)-2\big]dt-\int_{\check t}^{t_2}e^{-\lambda t}a(t-t_2)dt.
	    		 	\end{align*}
	    		 	Denote the last two integrals respectively by $I_1$ and $I_2$. Then, $J(\check x, \check u)|_{[t_1, t_2]}=I_1-I_2$. By regrouping and  applying the formula for integration by parts, one has
	    		 	\begin{align*}
	    		 	I_1& = -\dfrac{a}{\lambda}\int_{t_1}^{\check t}(t-t_1)d(e^{-\lambda t})-2\int_{t_1}^{\check t}e^{-\lambda t}dt \\
	    		 	& = -\dfrac{a}{\lambda}\big[(t-t_1)e^{-\lambda t}\big]\big|_{t_1}^{\check t}-\int_{t_1}^{\check t}e^{-\lambda t}dt\big]-2\int_{t_1}^{\check t}e^{-\lambda t}dt \\
	    		 	& = (\dfrac{a}{\lambda}-2)\int_{t_1}^{\check t}e^{-\lambda t}dt-\dfrac{a}{2\lambda}(t_2-t_1)e^{-\lambda \check t}\\
	    		 	& = (\dfrac{2}{\lambda}-\dfrac{a}{\lambda^2})\big(e^{-\lambda \check t}-e^{-\lambda t_1}\big)-\dfrac{a}{2\lambda}(t_2-t_1)e^{-\lambda \check t}.
	    		 	\end{align*}
	    		 	Similarly,
	    		 	\begin{align*}
	    		 	I_2 & = -\dfrac{a}{\lambda}\int_{\check t}^{t_2}(t-t_2)d(e^{-\lambda t})\\
	    		 	& = \dfrac{a}{\lambda}\big[(t-t_2)e^{-\lambda t}\big]\big|_{t_2}^{\check t}+\dfrac{a}{\lambda}\big[\int_{\check t}^{t_2}e^{-\lambda t}dt\big]\\
	    		 	& = -\dfrac{a}{2\lambda}(t_2-t_1)e^{-\lambda \check t}-\dfrac{a}{\lambda^2}\big[e^{-\lambda t_2}-e^{-\lambda \check t}\big].
	    		 	\end{align*}
	    		 	Thus, \begin{align*}
	    		 	J(\check x, \check u)|_{[t_1, t_2]} & = (\dfrac{2}{\lambda}-\dfrac{a}{\lambda^2})\big(e^{-\lambda \check t}-e^{-\lambda t_1}\big)-\dfrac{a}{2\lambda}(t_2-t_1)e^{-\lambda \check t}+\dfrac{a}{2\lambda}(t_2-t_1)e^{-\lambda \check
	    		 		t}\\
	    		 	& \quad +\dfrac{a}{\lambda^2}\big[e^{-\lambda t_2}-e^{-\lambda \check t}\big]\\
	    		 	& = (\dfrac{2}{\lambda}-\dfrac{2a}{\lambda^2})e^{-\lambda \check t}+(\dfrac{a}{\lambda^2}-\dfrac{2}{\lambda})e^{-\lambda t_1}+\dfrac{a}{\lambda^2}e^{-\lambda t_2}.
	    		 	\end{align*}
	    		 	Therefore,
	    		 	\begin{eqnarray}\label{identities}
	    		 	\begin{array}{rl}
	    		 	J(\check x, \check u)|_{[t_1, t_2]}-J(\widetilde x, \widetilde u)|_{[t_1, t_2]}
	    		 	& = \big(\dfrac{2}{\lambda}-\dfrac{2a}{\lambda^2}\big)e^{-\lambda \check t}+\big(\dfrac{a}{\lambda^2}-\dfrac{1}{\lambda}\big)e^{-\lambda t_1}+\big(\dfrac{a}{\lambda^2}-\dfrac{1}{\lambda}\big)e^{-\lambda t_2}\\
	    		 	& = \dfrac{1}{\lambda}\big(\dfrac{a}{\lambda}-1\big)\big(e^{-\lambda t_1}-2e^{-\lambda \check t}+e^{-\lambda t_2}\big).
	    		 	\end{array}
	    		 	\end{eqnarray} Thus, formula \eqref{identity} is proved. To obtain the second assertion of the lemma, put  $\psi(t)=e^{-\lambda t}$ for all $t\in\R$. Since $\psi''(t)>0$ for every~$t$, the function $\psi$ is strictly convex. So,
	    		 	$$\psi\big(\frac{1}{2}t_1+\frac{1}{2}t_2\big)<\frac{1}{2}\psi(t_1)+\frac{1}{2}\psi(t_2).$$
	    		 	It follows that $\Delta(t_1, t_2)>0$ for any $t_1<t_2$. Combining this with  \eqref{identities} and the inequality  $\dfrac{a}{\lambda}-1>0$, we obtain the strict inequality $J(\check x, \check u)|_{[t_1, t_2]}>J(\widetilde x, \widetilde u)|_{[t_1, t_2]}$.
	    		 \end{proof}
	    		 
	    		 The following analogue of Lemma \ref{lemma_technical1} will be used latter on. 
	    		 
	    		 \begin{Lemma}\label{lemma_technical1b}
	    		Let $t_1, t_2$ be as in Lemma \ref{lemma_technical1}. Let $J(x, u)|_{[t_1, t_2]}$ and $\Delta(t_1, t_2)$ be defined, respectively, by~\eqref{integral_on_segment} and \eqref{Delta}.	If $(\widetilde x, \widetilde u)$ and $(\hat x, \hat u)$ are feasible processes for $(FP_{3a})$ with $\widetilde x_1(t)=-1$ for all $t\in [t_1, t_2]$ and 
	    		 	\begin{equation*}\label{hat_x1}
	    		 	\hat x_1(t)=
	    		 	\begin{cases}
	    		 	-1+a(t-t_1), &\quad t\in [t_1, \hat t ]\\
	    		 	-1-a(t-t_2), &\quad t\in (\hat t, t_2],
	    		 	\end{cases}
	    		 	\end{equation*}
	    		 	where $\hat t:=2^{-1}(t_1+t_2)$, then one has 
	    		 	\begin{eqnarray*}\label{identity-1}
	    		 	J(\hat x, \hat u)|_{[t_1, t_2]}-J(\widetilde x, \widetilde u)|_{[t_1, t_2]}=-\dfrac{1}{\lambda}\big(\dfrac{a}{\lambda}-1\big)\Delta(t_1, t_2).
	    		 	\end{eqnarray*} Therefore,  $J( \hat x,  \hat u)|_{[t_1, t_2]}<J(\widetilde x, \widetilde u)|_{[t_1, t_2]}$.
	    		 \end{Lemma}
	    		 \begin{proof} By \eqref{state control system_FP_3a}, from our assumptions it follows that $\widetilde u(t)=0$ for almost all $t\in [t_1, t_2]$ and 
	    		 	\begin{equation*}
	    		 \hat u(t)=
	    		 	\begin{cases}
	    		 	-1, &\quad \mbox{a.e.}\ \, t\in [t_1,  \hat t ]\\
	    		 	1, &\quad \mbox{a.e.}\ \, t\in (\hat t, t_2].
	    		 	\end{cases}
	    		 	\end{equation*}
	    		 	Since $ \hat x (\cdot)$ is a feasible trajectory for $(FP_{3a})$, one has $ \hat x( \hat t)\leq 1$, i.e., $t_2-t_1\leq 4a^{-1}$. One has 
	    		 	\begin{equation*}
	    		 	J(\widetilde x, \widetilde u)|_{[t_1, t_2]}=\int_{t_1}^{t_2} \big[-e^{-\lambda t}\big(\widetilde x_1(t)+\widetilde u(t)\big)\big] dt=\int_{t_1}^{t_2} e^{-\lambda t}dt=-\dfrac{1}{\lambda}e^{-\lambda t_2}+\dfrac{1}{\lambda}e^{-\lambda t_1}.
	    		 	\end{equation*}
	    		 	Besides, the formulas for $\hat  x_1$ and $\hat  u$ on $[t_1, t_2]$ imply that
	    		 	\begin{align*}\
	    		 	J(\hat x, \hat u)|_{[t_1, t_2]}=& \int_{t_1}^{\hat  t}\big[-e^{-\lambda t}\big((-1+a(t-t_1))-1 \big)\big] dt\\& +\int_{\hat  t}^{t_2} \big[-e^{-\lambda t}\big((-1-a(t-t_2)) +1  \big)\big] dt\\
	    		 	=& -\int_{t_1}^{\hat  t}e^{-\lambda t}\big[a(t-t_1)-2\big]dt+\int_{\hat  t}^{t_2}e^{-\lambda t}a(t-t_2)dt.
	    		 	\end{align*}
	    		 	Thus, changing the sign of the expression $J(\hat x, \hat u)|_{[t_1, t_2]}-J(\widetilde x, \widetilde u)|_{[t_1, t_2]}$ we get the expression on the left-hand-side of \eqref{identity}. So, the desired results follow from Lemma~\ref{lemma_technical1}.
	    		 \end{proof}
	    		 
	    		 We will need two more lemmas.
	    		 
	    		 	\begin{Lemma}\label{lemma_technical2} Consider the function $\Delta:\R^2\to \R$ defined by \eqref{Delta}. For any
	    		 		$t_1, t_2\in\R$ with $t_1<t_2$ and for any $\bar\varepsilon \in (0, t_2-t_1)$, one has 
	    		 		\begin{equation}\label{Delta_1st_ineq}
	    		 		\Delta(t_1+\bar\varepsilon, t_2)<\Delta(t_1, t_2)
	    		 		\end{equation} and  
	    		 		\begin{equation}\label{Delta_2nd_ineq}
	    		 		\Delta(t_1, t_2)>\Delta(t_1, t_1+\bar\varepsilon)+\Delta(t_1+\bar\varepsilon, t_2).
	    		 		\end{equation}
	    		 		\end{Lemma}
	    		 	\begin{proof} Fix a value $\bar\varepsilon \in (0, t_2-t_1)$. To obtain \eqref{Delta_1st_ineq}, consider the function $\psi_1(\varepsilon):=\Delta(t_1+\varepsilon, t_2)$ of the variable $\varepsilon\in \R$. Since $\psi_1(\varepsilon)=e^{-\lambda (t_1+\varepsilon)}-2e^{-\frac{1}{2}\lambda (t_1+\varepsilon+t_2)}+e^{-\lambda t_2},$ one sees that
	    		 	$\psi_1(.)$ is continuously differentiable on $\R$ and 
	    		 	$\psi_1'(\varepsilon)=\lambda\big(e^{-\frac{1}{2}\lambda(t_1+\varepsilon+t_2)}-e^{-\lambda (t_1+\varepsilon)}\big).$   		 	
	    		 		As the function $r(t):= e^{-\lambda t}$ is strictly decreasing on $\R$, the last equality implies that  $\psi_1'(\varepsilon)<0$ for every $\varepsilon\in [0, t_2-t_1)$. Hence, the function  
	    		 		$\psi_1(.)$ is strictly decreasing on $[0, t_2-t_1)$. So, the inequality~\eqref{Delta_1st_ineq} is valid.
	    		 	
	    		 	To obtain \eqref{Delta_2nd_ineq}, observe from \eqref{Delta} that
	    			 \begin{align*}
	    			 & \Delta(t_1, t_2)-\Delta(t_1, t_1+\bar\varepsilon)-\Delta(t_1+\bar\varepsilon, t_2)\\
	    			 & = \ e^{-\lambda t_1}-2e^{-\lambda \frac{t_1+t_2}{2}}+e^{-\lambda t_2} - \big(e^{-\lambda t_1}-2e^{-\lambda (t_1+\frac{\bar\varepsilon}{2})}+e^{-\lambda (t_1+\bar\varepsilon)}\big)\\
	    			 & \quad - \big(e^{-\lambda (t_1+\bar\varepsilon)} -2e^{-\lambda (\frac{t_1+t_2}{2}+\frac{\bar\varepsilon}{2})}+e^{-\lambda t_2}\big)\\
	    			 & = \ 2\big[e^{-\lambda (\frac{t_1+t_2}{2}+\frac{\bar\varepsilon}{2})}- e^{-\lambda \frac{t_1+t_2}{2}}\big]-2\big[e^{-\lambda (t_1+\bar\varepsilon)}-e^{-\lambda (t_1+\frac{\bar\varepsilon}{2})}\big]
	    			 \end{align*}
	    	Applying the classical mean value theorem to the differentiable function $r(t)=e^{-\lambda t}$, one can find $\tau_1\in (t_1+\frac{\bar\varepsilon}{2}, t_1+\bar\varepsilon)$ and $\tau_2 \in (\frac{t_1+t_2}{2}, \frac{t_1+t_2}{2}+\frac{\bar\varepsilon}{2})$ such that
	    	\begin{equation*}
	    	e^{-\lambda (t_1+\bar\varepsilon)}-e^{-\lambda (t_1+\frac{\bar\varepsilon}{2})}=\frac{\bar\varepsilon}{2}(-\lambda)e^{-\lambda \tau_1},\end{equation*} \begin{equation*} e^{-\lambda (\frac{t_1+t_2}{2}+\frac{\bar\varepsilon}{2})}- e^{-\lambda \frac{t_1+t_2}{2}}= \frac{\bar\varepsilon}{2}(-\lambda)e^{-\lambda \tau_2}.
	    	\end{equation*}
	    	Thus, $\Delta(t_1, t_2)-\Delta(t_1, t_1+\bar\varepsilon)-\Delta(t_1+\bar\varepsilon, t_2) = \bar\varepsilon\lambda\big[e^{-\lambda \tau_1}-e^{-\lambda \tau_2}\big].$
	    	As the function $r(t)$ is strictly decreasing on $\R$ and $\tau_1<\tau_2$, one gets $e^{-\lambda \tau_1}-e^{-\lambda \tau_2}>0$; hence the inequality \eqref{Delta_2nd_ineq} is proved.
	    \end{proof}   
	     		
	    		 \begin{Lemma}\label{lemma_technical4}
	    		 Let there be given $t_1, t_2\in [t_0,T]$, $t_1<t_2$, and $\xi >0$. Suppose that $(\widetilde x^{\xi}, \widetilde u^{\xi})$ and $(\check x^\xi, \check u^\xi)$  are feasible processes for $(FP_{3a})$ with $\widetilde x_1^\xi(t)=\xi$ for all $t\in [t_1, t_2]$ and 
	    		 	    		 	\begin{equation}\label{hat-x1}
	    		 	    		 	\check x_1^\xi(t)=
	    		 	    		 	\begin{cases}
	    		 	    		 	\xi-a(t-t_1), &\quad t\in [t_1, \check t ]\\
	    		 	    		 	\xi+a(t-t_2), &\quad t\in (\check t, t_2],
	    		 	    		 	\end{cases}
	    		 	    		 	\end{equation}
	    		 	    		 	where $\check t:=2^{-1}(t_1+t_2)$. Then one has 
	    		 	    		 	\begin{eqnarray}\label{modified_identity}
	    		 	    		 		J(\check x^\xi, \check u^\xi)|_{[t_1, t_2]}-J(\widetilde x^\xi, \widetilde u^\xi)|_{[t_1, t_2]}=\dfrac{1}{\lambda}\big(\dfrac{a}{\lambda}-1\big)\Delta(t_1, t_2),
	    		 	    		 	\end{eqnarray} with $J(x, u)|_{[t_1, t_2]}$ and $\Delta(t_1, t_2)$ being defined respectively by \eqref{integral_on_segment} and \eqref{Delta}.  Besides, the strict inequality $J(\check x^\xi, \check u^\xi)|_{[t_1, t_2]}>J(\widetilde x^\xi, \widetilde u^\xi)|_{[t_1, t_2]}$ is valid. 
	    		 \end{Lemma}
	    		 \begin{proof} The proof is similar to that of Lemma~\ref{lemma_technical1}.   		
	    		 \end{proof}
	    		 
	    		       		 \begin{Proposition}\label{lemma-lower-boundary1}
	         	    		 The situation where $\bar x_1(t)=-1$ for all $t$ from a subsegment $[t_1, t_2]$ of $[t_0, T]$ with $t_1<t_2$ cannot happen.
	         	    		 \end{Proposition}
	         	    		\begin{proof}
	         	    		 Since $(\bar x, \bar u)$ is a $W^{1,1}$ local minimizer of $(FP_{3a})$, by Definition \ref{local_minimizer} there exists $\delta>0$ such that the process $(\bar x, \bar u)$ minimizes the quantity $g(x(t_0), x(T))=x_2(T)$ over all feasible processes $(x, u)$ of $(FP_{3a})$ with $\|\bar x-x\|_{W^{1,1}([t_0, T];\R^2)} \leq \delta$. 
	         	    		 
	         	    		 To prove our assertion, suppose on the contrary that there are $t_1, t_2$ with $t_0 \leq t_1 <  t_2 \leq T$ such that $\bar x_1(t)=-1$ for all $ t\in[t_1, t_2]$.  Fixing a number $\varepsilon\in (0, t_2- t_1)$, we consider the pair of functions $(\hat  x^\varepsilon, \hat  u^\varepsilon)$, where
	         	    		 	    	    		\begin{align*}
	         	    		 	    	    		\hat  x_1^\varepsilon (t):=
	         	    		 	    	    		\begin{cases}
	         	    		 	    	    		\bar x_1(t),  \ \;  & t\in [t_0, t_1) \cup (t_1+\varepsilon, T]\\
	         	    		 	    	    		-1+a(t-t_1), & t\in [t_1, t_1+2^{-1}\varepsilon]\\
	         	    		 	    	    		-1-a(t-t_1-\varepsilon), & t\in (t_1+2^{-1}\varepsilon, t_1+\varepsilon]
	         	    		 	    	    		\end{cases}
	         	    		 	    	    		\end{align*}
	         	    		and $\hat  u^\varepsilon(t):=-a^{-1}\dfrac{d\hat  x_1^\varepsilon (t)}{dt} $ for almost all $t\in [t_0,T]$. Clearly, $(\hat  x^\varepsilon, \hat  u^\varepsilon)$ is a feasible process of $(FP_{3a})$. By \eqref{phase_second_component_FP3a}, \eqref{integral_on_segment}, and the definition of $\hat  x_1^\varepsilon (.)$,  we have 
	         	    	   \begin{align}\label{minus-lower-boundary} \bar x_2(T)-\hat  x_2^\varepsilon (T) = J(\bar x, \bar u)_{|[t_1, t_1+\varepsilon]}-J(\hat  x^\varepsilon, \hat  u^\varepsilon)_{|[t_1, t_1+\varepsilon]}.
	         	    		\end{align}
	         	    	Besides, it follows from Lemma~\ref{lemma_technical1b} and the constructions of $\bar x$ and $\hat  x^\varepsilon$ on $[t_1, t_1+\varepsilon]$ that
	         	    		\begin{align*}
	         	    		J(\bar x, \bar u)_{|[t_1, t_1+\varepsilon]}-J(\hat  x^\varepsilon, \hat  u^\varepsilon)_{|[t_1, t_1+\varepsilon]}>0.
	         	    		\end{align*}
	         	    		Combining this with \eqref{minus-lower-boundary} yields $\bar x_2(T)>\hat  x_2^\varepsilon (T)$, which contradicts  the $W^{1,1}$ local optimality of $(\bar x, \bar u)$, because $\|\bar x-\hat  x^\varepsilon\|_{W^{1,1}([t_0, T];\R^2)}\leq\delta$ for $\varepsilon >0$ small enough.  
	         	    		 \end{proof}
	         	    		
	    		  The following two propositions are crucial for describing the behavior of the local solutions of  $(FP_{3a})$.
	    		    		    	
	    	\begin{Proposition}\label{lemma-lower-boundary2}
	    	One must have $\bar x_1(t)>-1$ for all $ t \in (t_0, T)$. 
	    		\end{Proposition}
	    		\begin{proof} By our standing assumption, $(\bar x, \bar u)$ is a $W^{1,1}$ a local minimizer for $(FP_{3a})$. Let $\delta>0$ be chosen as in the proof of Proposition~\ref{lemma-lower-boundary1}. If the assertion is false, there would exist $\check t \in (t_0, T)$ with $\bar x_1(\check t)=-1$. 
	    			
	    			If there are $\varepsilon_1>0$ and $\varepsilon_2>0$ such that  $\bar x_1(t)>-1$ for all $t\in (\check t-\varepsilon_1,\check t)\cup (\check t,\check t+\varepsilon_2)$. Then, thanks to the continuity of $\bar x_1(.)$, by shrinking $\varepsilon_1>0$ and $\varepsilon_2>0$ (if necessary) one may assume that  $\bar x_1(t)\in (-1,1)$ for all $t\in (\check t-\varepsilon_1,\check t)\cup (\check t,\check t+\varepsilon_2)$.
	    		Then, since the curve $\bar x_1(.)$ cannot have more than one turning on the interval $(\check t-\varepsilon_1,\check t)$ (resp., on the interval $(\check t,\check t+\varepsilon_2)$) by the observation given at the beginning of the proof of Proposition~\ref{Lemma_main_tool}. So, replacing $\varepsilon_1$  (resp., $\varepsilon_2$) by a smaller positive number, one may assume that
	    		\begin{equation}\label{curve_down-up}
	    		\bar x_1(t)=
	    		\begin{cases}
	    		-1-a(t-\check t), &\quad t\in [\check t-\varepsilon_1, \check t ]\\
	    		-1+a(t-\check t), &\quad t\in (\check t,\check t+\varepsilon_2].
	    		\end{cases}
	    		\end{equation} To get a contradiction, we can apply the construction given in Lemma~\ref{lemma_technical4}. Namely, choose $\varepsilon >0$ as small as $\varepsilon<\min\{ \varepsilon_1,\varepsilon_2\}$ 
	    		and define a feasible process $(\widetilde x^{\varepsilon}, \widetilde u^{\varepsilon})$ for $(FP_{3a})$ by setting 
	    		\begin{equation}\label{tilde-u-epsilon}
	    	\widetilde u^{\varepsilon}(t)=
	    		\begin{cases}
	    		0, &\quad t\in [\check t-\varepsilon, \check t+\varepsilon]\\
	    		\bar u(t), &\quad t\in [t_0,\check t-\varepsilon)\cup (\check t+\varepsilon,T]
	    		\end{cases}
	    		\end{equation}
	    		and
	    			\begin{equation}\label{tilde-x1-epsilon}
	    			\widetilde x^{\varepsilon}(t)=
	    			\begin{cases}
	    			\bar x_1(\check t-\varepsilon), &\quad t\in [\check t-\varepsilon, \check t+\varepsilon]\\
	    			\bar x(t), &\quad t\in [t_0,\check t-\varepsilon)\cup (\check t+\varepsilon,T].
	    			\end{cases}
	    			\end{equation}
	    	Then, by  Lemma~\ref{lemma_technical4} one has 
	    		 $J(\bar x, \bar u)>J(\widetilde x^\varepsilon, \widetilde u^\varepsilon)$. This contradicts the $W^{1,1}$ local optimality of $(\bar x, \bar u)$, because $\|\bar x-\hat  x^\varepsilon\|_{W^{1,1}([t_0, T];\R^2)}\leq\delta$ for $\varepsilon >0$ small enough.  
	    		 
	    		 Since one cannot find $\varepsilon_1>0$ and $\varepsilon_2>0$ such that the strict inequality $\bar x_1(t)>-1$ holds for all $t\in (\check t-\varepsilon_1,\check t)\cup (\check t,\check t+\varepsilon_2)$, there must exist a sequence $\{t_k\}$ in $(t_0,T)$ converging to $\check t$ such that either $t_k<\check t$ for all $k$ or $t_k>\check t$ for all $k$, and $\bar x_1(t_k)=-1$ for each $k$. It suffices to consider the case $t_k<\check t$ for all $k$, as the other case can be treated similarly. By considering a subsequence (if necessary), we may assume that $t_k<t_{k+1}$ for all $k$.
	    		 
	    		 Choose $\bar k$ as large as 
	    		 \begin{equation}\label{the-height}
	    		 \check t-t_{\bar k}<\min \{2\delta a^{-1},\, 4 a^{-1}\}.
	    		 	\end{equation} 
	    		 	This choice of $\bar k$ guarantees that $\bar x_1(t)<1$ for every $t\in [t_{\bar k},\check t]$. Indeed, otherwise there is some $\alpha\in  (t_{\bar k},\check t)$ with $\bar x_1(\alpha)=1$. 
	    		  Setting $$\alpha_1=\min\big\{t\in [t_{\bar k},\alpha]\,:\, \bar x_1(t)=1\big\},\  \,\alpha_2=\max\big\{t\in [\alpha,\check t]\,:\, \bar x_1(t)=1\big\},$$ one has $\alpha_1\leq\alpha_2$, $[\alpha_1,\alpha_2]\subset [t_{\bar k},\check t]$, and $\bar x_1(t)\in (-1,1)$ for all $t\in (t_{\bar k},\alpha_1)\cup (\alpha_2,\check t) $. Then, by assertion S1 of Proposition~\ref{Lemma_main_tool}, one has $\alpha_1-t_{\bar k}=2a^{-1}$. Similarly, by assertion~S2 in that proposition, one has
	    		 	$\check t-\alpha_2=2a^{-1}$. So, one gets $\check t-t_{\bar k}\geq 4a^{-1}$, which comes in conflict with~\eqref{the-height}. 
	    		 	
	    		 	By Proposition~\ref{lemma-lower-boundary1}, one cannot have $\bar x_1(t)=-1$ for all $t\in [t_{\bar k},t_{\bar k+1}]$. Thus, there is some $\tau\in (t_{\bar k},t_{\bar k+1})$ with $\bar x_1(\tau)>-1$. Setting $$\tau_1=\max\big\{t\in [t_{\bar k},\tau]\,:\, \bar x_1(t)=-1\big\},\  \,\tau_2=\min\big\{t\in [\tau,t_{\bar k+1}]\,:\, \bar x_1(t)=-1\big\},$$ one has $\tau_1<\tau_2$, $[\tau_1,\tau_2]\subset [t_{\bar k},t_{\bar k+1}]$, and $\bar x_1(t)\in (-1,1)$ for all $t\in (\tau_1,\tau_2)$. Hence, replacing $t_{\bar k}$ (resp., $t_{\bar k+1}$) by $\tau_1$ (resp., $\tau_2$), one sees that all the above-described properties of the sequence $\{t_k\}$ remain and, in addition,
	    		 	 \begin{equation}\label{first_interval}
	    		 	 \bar x_1(t)\in (-1,1),\quad \forall t\in (t_{\bar k},t_{\bar k+1}).
	    		 	 \end{equation}
	    		 		    		 	
	    		 	Let $F:=\{t\in [t_{\bar k},\check t]\,:\, \bar x_1(t)=-1\}$ and $E:= [t_{\bar k},\check t]\setminus F$. Since $F$ is a closed subset of $\R$ and $E= (t_{\bar k},\check t)\setminus F$, $E$ is an open subset of $\R$. So, $E$ is the union of a countable family of disjoint open intervals (see \cite[Proposition~9, p.~17]{Royden_Fitzpatrick_2010}). Since $t_k\notin E$ for all $k$, we have a representation $E=\displaystyle\bigcup_{i=1}^\infty E_i$, where the intervals $E_i=(\tau_1^{(i)},\tau_2^{(i)}),\ i\in\N$, are nonempty and disjoint. Thanks to \eqref{first_interval}, one may suppose that $E_1=(\tau_1^{(1)},\tau_2^{(1)})=(t_{\bar k},t_{\bar k+1})$.  Note also that, for any $i\in\N$, $ \bar x_1(t)\in (-1,1)$ for all $t\in E_i$.  Since $\bar x_1(\tau_1^{(i)})=\bar x_1(\tau_2^{(i)})=-1$, by assertion S3 of Proposition~\ref{Lemma_main_tool} one gets
	    		 	\begin{equation}\label{bar_x1_Ei}
	    		 	\bar x_1(t)=
	    		 	\begin{cases}
	    		 	-1+a(t-\tau_1^{(i)}), \quad & t\in [\tau_1^{(i)}, 2^{-1}(\tau_1^{(i)}+\tau_2^{(i)})]\\
	    		 -1-a(t-\tau_2^{(i)}), & t\in (2^{-1}(\tau_1^{(i)}+\tau_2^{(i)}), \tau_2^{(i)}].
	    		 	\end{cases}
	    		 	\end{equation}
	    		 	    		 	
	    		 	If  the set $F_1:=F\setminus \{t_{\bar k}\}$ has an isolated point in the induced topology of $[t_{\bar k},\check t]$, says, $\bar t$. Then, one must have  $\bar t\in [t_{\bar k+1},\check t)$. So, there exists $\varepsilon>0$ such that $(\bar t-\varepsilon,\bar t+\varepsilon)\subset (t_{\bar k},\check t)$ and $\bar x_1(t)\in (-1,1)$ for all $t\in (\bar t-\varepsilon,\bar t)\cup (\bar t,\bar t+\varepsilon)$. Applying the construction given in the first part of this proof, we find a feasible process $(\widetilde x^{\varepsilon}, \widetilde u^{\varepsilon})$ for $(FP_{3a})$ with the property $J(\bar x, \bar u)>J(\widetilde x^\varepsilon, \widetilde u^\varepsilon)$. This contradicts the $W^{1,1}$ local optimality of $(\bar x, \bar u)$, because \eqref{the-height} assures that $\|\bar x-\hat  x^\varepsilon\|_{W^{1,1}([t_0, T];\R^2)}\leq\delta$.  
	    		 	
	    		 	Now, suppose that every point in the compact set $F_1$ is a limit point of this set in the induced topology of $[t_{\bar k},\check t]$. Then, if the Lebesgue measure $\mu_L(F_1)$ of $F_1$ is null, then the structure of $F_1$ is similar to that of \textit{the Cantor set}\footnote{https://en.wikipedia.org/wiki/Cantor$\_$set.}, constructed from the segment $[t_{\bar k+1},\check t]\subset\R$. If $\mu_L(F_1)>0$, the structure of $F_1$ is similar to that of a \textit{fat Cantor set}, which is also called \textit{a Smith-Volterra-Cantor set}\footnote{https://en.wikipedia.org/wiki/Smith-Volterra-Cantor$\_$set.}.
	    		 	
	    		 	Putting
	    		 	\begin{equation}\label{tilde-u-lower}
	    		 	\widetilde u(t)=
	    		 	\begin{cases}
	    		 	0, &\quad t\in [t_{\bar k},\check t]\\
	    		 	\bar u(t), &\quad t\in [t_0,t_{\bar k})\cup (\check t,T]
	    		 	\end{cases}
	    		 	\end{equation}
	    		 	and
	    		 	\begin{equation}\label{tilde-x1-lower}
	    		 	\widetilde x_1(t)=
	    		 	\begin{cases}
	    		 	-1, & \quad t\in  [t_{\bar k},\check t]\\
	    		 	\bar x_1(t), &\quad t\in  [t_0,t_{\bar k})\cup (\check t,T],
	    		 	\end{cases}
	    		 	\end{equation} we see that $(\widetilde x,\widetilde u)$ is a feasible process for $(FP_{3a})$. Similarly, define
	    		 	\begin{equation}\label{Big-hat-u-1}
	    		 	u(t)=
	    		 	\begin{cases}
	    		 	-1, &\quad t\in [t_{\bar k+1},2^{-1}(t_{\bar k+1}+\check t)]\\
	    		 	1, &\quad t\in (2^{-1}(t_{\bar k+1}+\check t),\check t]\\
	    		 	\bar u(t), &\quad t\in [t_0,t_{\bar k+1})\cup (\check t,T]
	    		 	\end{cases}
	    		 	\end{equation}
	    		 	and
	    		 	\begin{equation}\label{Big-hat-x1-1}
	    		 	x_1(t)=
	    		 	\begin{cases}
	    		 	-1+a(t-t_{\bar k+1}), &\quad t\in [t_{\bar k+1},2^{-1}(t_{\bar k}+\check t)]\\
	    		 	-1-a(t-\check t), &\quad t\in (2^{-1}(t_{\bar k+1}+\check t),\check t]\\
	    		 	\bar x_1(t), &\quad t\in [t_0,t_{\bar k+1})\cup (\check t,T],
	    		 	\end{cases}
	    		 	\end{equation} and observe that $(x,u)$ is a feasible process for $(FP_{3a})$. Using \eqref{the-height}, it is easy to verify that $\|x-\bar x\|_{W^{1,1}([t_0, T];\R^2)} \leq \delta$. Thus, if it can be shown that
	    		 	\begin{equation}\label{difference_Jhatx_Jbarx}
	    		 	J(x,u)<J(\bar x,\bar u),
	    		 	\end{equation} then we get a contradiction to the $W^{1,1}$ local optimality of $(\bar x, \bar u)$. Hence, the proof of the lemma will be completed.
	    		 	
	    		 	By \eqref{tilde-u-lower}--\eqref{Big-hat-x1-1} and Lemma \ref{lemma_technical1b}, one has $J(\widetilde x, \widetilde u)-J(x,u) = J(\widetilde x, \widetilde u)|_{[t_{\bar k},\check t]}-J(x,u)|_{[t_{\bar k},\check t]}$. Therefore,
	    		 		\begin{eqnarray}\label{identity_Jhatx_Jtildex}
	    		 		 J(\widetilde x, \widetilde u)-J(x,u) =\dfrac{1}{\lambda}\big(\dfrac{a}{\lambda}-1\big)\big[
	    		 			\Delta(t_{\bar k},t_{\bar k+1})+\Delta(t_{\bar k+1},\check t)\big],
	    		 		\end{eqnarray} where $\Delta(t_1, t_2)$, for any $t_1, t_2$ with  $t_1<t_2$,  is given by \eqref{Delta}. In addition, using  \eqref{tilde-u-lower}, \eqref{tilde-x1-lower}, the decomposition $[t_{\bar k+1},\check t]=\big(\displaystyle\bigcup_{i=2}^\infty E_i\big)\cup F_1$, and the sum rule \cite[Theorem~1', p. 297]{Kolmogorov_Fomin_1970} and the decomposition formula \cite[Theorem~4, p. 298]{Kolmogorov_Fomin_1970}  for the Lebesgue integrals, one gets
	    		 			\begin{align*}
	    		 			J(\bar x,\bar u)-J(\widetilde x, \widetilde u)  & = J(\bar x,\bar u)|_{[t_{\bar k},\check t]}-J(\widetilde x, \widetilde u)|_{[t_{\bar k},\check t]}\nonumber\\
	    		 			& = \displaystyle\int_{[t_{\bar k},\check t]}\Big[-e^{-\lambda t}(\big[\bar x_1(t)+\bar u(t)\big]-\big[\widetilde x_1(t)+\widetilde u(t)]\big)\Big]dt\\
	    		 			& = \sum_{i=2}^{\infty} \displaystyle\int_{E_i}\Big[-e^{-\lambda t}(\big[\bar x_1(t)+\bar u(t)\big]-\big[\widetilde x_1(t)+\widetilde u(t)]\big)\Big]dt\\
	    		 			&  + \displaystyle\int_{F_1}\Big[-e^{-\lambda t}(\big[\bar x_1(t)+\bar u(t)\big]-\big[\widetilde x_1(t)+\widetilde u(t)]\big)\Big]dt.	\end{align*} Hence, it holds that
	    		 			\begin{align}\label{identity_Jbarx_Jtildex}
	    		 			J(\bar x,\bar u)-J(\widetilde x, \widetilde u)  = - \dfrac{1}{\lambda}\big(\dfrac{a}{\lambda}-1\big) \sum_{i=2}^{\infty}\Delta\big(\tau_1^{(i)},\tau_2^{(i)}\big) +I,
	    		 		 \end{align}
	    		 		where $I:= \displaystyle\int_{F_1}\Big[-e^{-\lambda t}(\big[\bar x_1(t)+\bar u(t)\big]-\big[\widetilde x_1(t)+\widetilde u(t)]\big)\Big]dt.$	 Given any $t\in F_1$, we observe that $\bar x_1(t)=\widetilde x_1(t)=-1$ and $\widetilde u(t)=0$. Since every point in $F_1$ is a limit point of this set in the induced topology of $[t_{\bar k},\check t]$, we can find a sequence $\{\xi_j^t\}$ in $F_1$ satisfying $\displaystyle\lim_{j\to\infty}\xi_j^t=t$. As the derivative $\bar x_1(t)$ exists a.e. on $[t_0,T]$, it exists a.e. on $F_1$. In combination with the first differential equation in \eqref{state control system_FP_3a}, this yields $\dot{\bar x}_1(t)=-a\bar u(t)$ a.e. $t\in F_1$. Since $\bar x_1(t)=-1$ for all $t\in F_1$, for a.e. $t\in F_1$ it holds that
	    		 		$$\bar u(t)=-\frac{1}{a}\,\dot{\bar x}_1(t)=-\frac{1}{a}\,\lim_{j\to\infty}\frac{{\bar x}_1(\xi_j^t)-{\bar x}_1(t)}{\xi_j^t-t}=0.$$ We have thus shown that $\big[\bar x_1(t)+\bar u(t)\big]-\big[\widetilde x_1(t)+\widetilde u(t)]=0$ for a.e. $t\in F_1$. This implies that $I=0$. Now, adding \eqref{identity_Jhatx_Jtildex} \eqref{identity_Jbarx_Jtildex}, we get  
	    		 		\begin{eqnarray}\label{identity_Jbarx_Jx}
	    		 		J(\bar x,\bar u)-J(x,u) =\dfrac{1}{\lambda}\big(\dfrac{a}{\lambda}-1\big)\Big[
	    		 		\Delta(t_{\bar k},t_{\bar k+1})+\Delta(t_{\bar k+1},\check t)-\sum_{i=2}^{\infty}\Delta\big(\tau_1^{(i)},\tau_2^{(i)}\big)\Big].
	    		 		\end{eqnarray} We have
	    		 		\begin{eqnarray}\label{inequality_Deltas}
	    		 		\sum_{i=2}^{\infty}\Delta\big(\tau_1^{(i)},\tau_2^{(i)}\big)
	    		 		\leq \Delta(t_{\bar k+1},\check t).
	    		 		\end{eqnarray} To establish this inequality, we first show that
	    		 		\begin{eqnarray}\label{inequality_partial_sums}
	    		 		\sum_{i=2}^{m}\Delta\big(\tau_1^{(i)},\tau_2^{(i)}\big)
	    		 		< \Delta(t_{\bar k+1},\check t)
	    		 		\end{eqnarray} for any integer $m\geq 2$. Taking account of the fact that every point in $F_1$ is a limit point of this set in the induced topology of $[t_{\bar k},\check t]$, by reordering the intervals $\big(\tau_1^{(i)},\tau_2^{(i)}\big)$ for $i=2,\dots,m$, we may assume that 
	    		 		$t_{\bar k+1}<\tau_1^{(2)}<\tau_2^{(2)}<\tau_1^{(3)}<\tau_2^{(3)}<\dots<\tau_1^{(m)}<\tau_2^{(m)}<\check t.$ Then, by Lemma~\ref{lemma_technical2} and by induction, we have
	    		 		\begin{eqnarray*}
	    		 			\sum_{i=2}^{m}\Delta\big(\tau_1^{(i)},\tau_2^{(i)}\big)&<&  \Big[\Delta\big(t_{\bar k+1},\tau_1^{(2)}) +		 			
	    		 			\Delta\big(\tau_1^{(2)},\tau_2^{(2)}\big)\Big]+	\sum_{i=3}^{m}\Delta\big(\tau_1^{(i)},\tau_2^{(i)}\big)\\
	    		 			& < &\Big[\Delta\big(t_{\bar k+1},\tau_2^{(2)})+\Delta\big(\tau_2^{(2)},\tau_1^{(3)}\big)\Big]+\sum_{i=3}^{m}\Delta\big(\tau_1^{(i)},\tau_2^{(i)}\big)\\
	    		 			& \vdots & \\
	    		 			& <&\Delta\big(t_{\bar k+1},\tau_2^{(m)})+\Delta\big(\tau_2^{(m)},\check t\big)\\
	    		 			& <& \Delta(t_{\bar k+1},\check t).
	    		 	  \end{eqnarray*} Thus, \eqref{inequality_partial_sums} is valid. Since $\Delta\big(\tau_1^{(i)},\tau_2^{(i)}\big)>0$ for all $i=2,3,\dots$, the estimate \eqref{inequality_partial_sums} shows that the series  	$\displaystyle\sum_{i=2}^{\infty}\Delta\big(\tau_1^{(i)},\tau_2^{(i)}\big)$ is convergent. Letting $m\to\infty$, from  \eqref{inequality_partial_sums} one obtains~\eqref{inequality_Deltas}. Since $\Delta(t_{\bar k},t_{\bar k+1})>0$, the equality \eqref{identity_Jbarx_Jx} and the inequality \eqref{inequality_Deltas} imply \eqref{difference_Jhatx_Jbarx}.  
	    		 	  
	    		 	  The proof is complete.   		 		
	    		\end{proof}
	    			
To continue, using the data set $\{a,\lambda, t_0,T,x_0\}$ of $(FP_{3a})$, we define $\rho=\dfrac{1}{\lambda}\ln \dfrac{a}{a-\lambda}>0$ and $\bar t=T-\rho$. Besides, for a given $x_0\in [-1, 1]$, let 
\begin{equation}\label{two_rhos}
\rho_1:=a^{-1}(1+x_0)\quad {\rm and}\quad  \rho_2:=a^{-1}(1-x_0).
\end{equation}
 As $x_0\in [-1,1]$, one has $\rho_1\in [0, 2a^{-1}]$ and $\rho_2\in [0, 2a^{-1}]$. Moreover, since $\bar x_1(t)$ is a continuous function,
	${\mathcal T}_1:=\{t\in [t_0,T]\, :\, \bar x_1(t)=1\}$ is a compact set (which may be empty). If ${\mathcal T}_1$ is nonempty, then we consider the numbers $\alpha_1:=\min \{t\,:\,t\in {\mathcal T}_1\}$ and $\alpha_2:=\max \{t\,:\,t\in {\mathcal T}_1\}$. 
	
	\medskip
	By Proposition~\ref{lemma-lower-boundary2}, one of next four cases must occur.
	  
	  \medskip
	  {\bf Case 1:} \textit{$\bar x_1(t)> -1$ for all $t\in [t_0, T]$}. Then, condition (i) means that \eqref{1st_condition_for_mu_FP3a} and~\eqref{2nd_condition_for_mu_FP3a} are satisfied, while conditions (ii)--(iv) remain the same as those in Subsection~4.2 of Part 1 (\cite{VTHuong_Yao_Yen_Part1}). So, the curve $\bar x_1(t)$ must have of one of the forms ${\rm (a)}$--${\rm (c)}$ depicted in Theorem~4.4 of Part 1 (\cite{VTHuong_Yao_Yen_Part1}), where we let $\bar x_1(t)$ play the role of $\bar x(t)$. Of course, the condition $\bar x_1(t)> -1$ for all $t\in [t_0, T]$ must be satisfied. Note that the latter is equivalent to the requirement $\bar x_1(T)> -1$. With respect to the just mentioned three forms of $\bar x(t)$, we have the following three subcases.
	  
	   \underline{\textit{Subcase 1a}}: $\bar x_1(t)$ is given by
	  \begin{equation}\label{traj-subcase1a-FP3a}
	  \bar x_1(t)=x_0-a(t-t_0), \quad t \in [t_0, T].
	  \end{equation}
	  By statement (a) of Theorem~4.4 of Part 1 (\cite{VTHuong_Yao_Yen_Part1}), this situation happens when $T-t_0\leq \rho$. By \eqref{traj-subcase1a-FP3a}, condition $\bar x_1(T)> -1$ is equivalent to $T-t_0<\rho_1$. Therefore, if either $\rho< \rho_1$ and $T-t_0\leq \rho$, or $\rho\geq \rho_1$ and $T-t_0< \rho_1$, then $\bar x_1(t)$ is given by \eqref{traj-subcase1a-FP3a}.
	  	
	  \underline{\textit{Subcase 1b}}: $\bar x_1(t)$ is given by
	  	    		 \begin{equation}\label{traj-subcase1b-FP3a}
	  	    		 \bar x_1(t)=
	  	    		 \begin{cases}
	  	    		 x_0+a(t-t_0), \quad & t \in [t_0, \bar t]\\
	  	    		 x_0-a(t+t_0-2\bar t), & t \in (\bar t, T].
	  	    		 \end{cases}
	  	    		 \end{equation}
	 Then, statement (b) of Theorem~4.4 of Part 1 (\cite{VTHuong_Yao_Yen_Part1}) requires that $\rho< T-t_0<\rho+\rho_2$. By \eqref{traj-subcase1b-FP3a}, the inequality $\bar x_1(T)>-1$ means $T-t_0>2\rho-\rho_1$. Thus, if $\max \{\rho; 2\rho-\rho_1\}< T-t_0<\rho+\rho_2$, then $\bar x_1(t)$ is given by \eqref{traj-subcase1b-FP3a}.
	  
\underline{\textit{Subcase 1c}}: $\bar x_1(t)$ is given by	    		 
	 
\begin{equation}\label{traj-subcase1c-FP3a}
\bar x_1(t)=
\begin{cases}
x_0+a(t-t_0), \quad & t \in [t_0, t_0+\rho_2]\\
1-a(t-t_0-\rho_2), & t \in (t_0+\rho_2, T].
\end{cases}
\end{equation}
Since $\alpha_1=t_0+a^{-1}(1-x_0)=t_0+\rho_2$, this situation is in full agreement with the one in assertion (c) of Theorem~4.4 of Part 1 (\cite{VTHuong_Yao_Yen_Part1}). Here, one must have $T-t_0\geq\rho+\rho_2$. By \eqref{traj-subcase1c-FP3a}, the inequality $\bar x_1(T)>-1$ means $T-t_0<a^{-1}(3-x_0)$. Thus, this situation occurs if $\rho+\rho_2\leq T-t_0<a^{-1}(3-x_0)$.
	  	
 {\bf Case 2:} \textit{$\bar x_1(t_0)=-1$ and $\bar x_1(t)>-1$ for all $t\in (t_0, T]$}. Let $\bar\varepsilon>0$ be such that $t_0+\bar\varepsilon< T$. For any $k\in\N$ with $k^{-1}\in (0,\bar\varepsilon)$, by the comments before Propositions~4.1 and by Proposition~4.2 of Part~1 (\cite{VTHuong_Yao_Yen_Part1}) we can assert that the restriction of $(\bar x, \bar u)$ on $[t_0+k^{-1}, T]$ is a $W^{1,1}$ local minimizer for the Mayer problem obtained from $(FP_{3a})$ by replacing $t_0$ with $t_0+k^{-1}$. Since $\bar x_1(t)>-1$ for all $t\in [t_0+k^{-1},T]$, repeating the arguments already used in Case~1 yields a formula for $\bar x_1(t)$ on $[t_0+k^{-1},T]$. With $\rho_1 (k):=a^{-1}[1+\bar x_1(t_0+k^{-1})]$ and  $\rho_2(k):=a^{-1}[1-\bar x_1(t_0+k^{-1})]$, for every $k\in\N$ we see that the function $\bar x_1(t)$ on $[t_0+k^{-1},T]$ must belong to one of the following three categories, which correspond to the three forms of the function $\bar x_1(t)$ in Case~1.
\begin{enumerate}[(C1)]
\item $\bar x_1(t)$ is given by $$\bar x_1(t)=\bar x_1(t_0+k^{-1})-a(t-t_0-k^{-1}), \quad   t\in [t_0+k^{-1}, T],$$
provided that $\rho< \rho_1(k)$ and $T-t_0-k^{-1}\leq \rho$, or $\rho\geq \rho_1(k)$ and $T-t_0-k^{-1}< \rho_1(k)$.
\item $\bar x_1(t)$ is given by
\begin{equation*}
\bar x_1(t)=
\begin{cases}
\bar x_1(t_0+k^{-1})+a(t-t_0-k^{-1}), \quad & t \in [t_0, \bar t]\\
\bar x_1(t_0+k^{-1})-a(t+t_0+k^{-1}-2\bar t), & t \in (\bar t, T],
\end{cases}
\end{equation*}
provided that $\max \{\rho; 2\rho-\rho_1(k)\}< T-t_0-k^{-1}<\rho+\rho_2(k)$.
\item $\bar x_1(t)$ is given by
\begin{equation*}
\bar x_1(t)=
\begin{cases}
\bar x_1(t_0+k^{-1})+a(t-t_0-k^{-1}), \quad & t \in [t_0, t_0+\rho_2(k)]\\
1-a(t-t_0-k^{-1}-\rho_2(k)), & t \in (t_0+\rho_2(k), T],
\end{cases}
\end{equation*}
provided that $\rho+\rho_2(k)\leq T-t_0-k^{-1}<a^{-1}[3-\bar x_1(t_0+k^{-1})]$.
\end{enumerate}	 

 By the Dirichlet principle, there exist an infinite number of indexes $k$ with $k^{-1}\in (0,\bar\varepsilon)$ such that the formula for $\bar x_1(t)$ is given in the category C1 (resp., C2, or C3). By considering a subsequence if necessary, we may assume that this happens for all $k$ with $k^{-1}\in (0,\bar\varepsilon)$.
 
If the first situation occurs, then by letting $k\to\infty$ we have $\bar x_1(t)=-1-a(t-t_0)$ for all $t\in [t_0, T]$. This is impossible since the requirement $\bar x_1(t)>-1$ for all $t\in (t_0, T]$ is violated.
 
If the second situation occurs, then by letting $k\to\infty$ we have 
	  	 	  	    		 \begin{equation}\label{traj-subcase2b-FP3a}
	  	 	  	    		 \bar x_1(t)=
	  	 	  	    		 \begin{cases}
	  	 	  	    		 -1+a(t-t_0), \quad & t \in [t_0, \bar t]\\
	  	 	  	    		 -1-a(t+t_0-2\bar t), & t \in (\bar t, T],
	  	 	  	    		 \end{cases}
	  	 	  	    		 \end{equation}
	  	 	 provided that $2\rho\leq T-t_0\leq\rho+2a^{-1}$. Since $\bar x_1(t)>-1$ for all $t\in (t_0, T]$, especially $\bar x_1(T)>-1$, one must have $2\rho < T-t_0$.

If the last situation occurs, then $\bar x_1(t)$ is given by	    		 
	  	  	    		 \begin{equation}\label{traj-subcase2c-FP3a}
	  	  	    		 \bar x_1(t)=
	  	  	    		 \begin{cases}
	  	  	    		  -1+a(t-t_0), \quad & t \in [t_0, t_0+2a^{-1}]\\
	  	  	    		 1-a(t-t_0-2a^{-1}), & t \in (t_0+2a^{-1}, T],
	  	  	    		 \end{cases}
	  	  	    		 \end{equation}
	  	  	provided that $\rho+2a^{-1}\leq T-t_0\leq 4a^{-1}$. Having in mind that $\bar x_1(T)>-1$, one must have the strict inequality $T-t_0< 4a^{-1}$.
	  	  	
Since the first situation cannot happen and since $\bar t=t_0+2a^{-1}$ when $T-t_0=\rho+2a^{-1}$, our results in this case can be summarized as follows.

\underline{\textit{Subcase~2a}}: $\bar x_1(t)$ is given by \eqref{traj-subcase2b-FP3a}, provided that $2\rho< T-t_0<\rho+2a^{-1}$.

\underline{\textit{Subcase 2b}}: $\bar x_1(t)$ is given by \eqref{traj-subcase2c-FP3a}, provided that $\rho+2a^{-1}\leq T-t_0<4a^{-1}$.
	  	
	  \smallskip   	 	   	  
	  {\bf Case 3:} \textit{$\bar x_1(T)=-1$ and $\bar x_1(t)>-1$ for all $t\in [t_0, T)$}. We split this case into two subcases.
	  
	 \underline{\textit{Subcase 3a}}: ${\mathcal T}_1=\emptyset$. Then $\bar x_1(t)\in (-1, 1)$ for all $t\in [t_0, T)$ and $\bar x_1(T)=-1$. By some arguments similar to those of the proof of Proposition~\ref{Lemma_main_tool}, one can show that formula for $\bar x_1(.)$ on $[t_0, T]$ is one of the following two types: 
	 \begin{equation}\label{type1}
	 \bar x_1(t)=x_0-a(t-t_0), \quad  t\in [t_0, T],
	 \end{equation}
	  and
	 \begin{equation}\label{type2}
	 	  \bar x_1(t)=
	 	  \begin{cases}
	 	  x_0+a(t-t_0), \quad & t\in [t_0, t_{\zeta}]\\
	 	  -1-a(t-T), &  t\in (t_{\zeta}, T],
	 	  \end{cases}
	 	  \end{equation}
	with $t_{\zeta}\in (t_0, T)$. 
	
	If $\bar x_1(.)$ is given by \eqref{type1}, then  $\bar x_1(T)=-1$ if and only if $T-t_0=\rho_1$. Since $x_0\in (-1, 1]$, the latter yields $0<T-t_0=\rho_1 \leq 2a^{-1}$. 
	
	If $\bar x_1(.)$ is of the form~\eqref{type2}, then the equality $\bar x_1(T)=-1$ implies that $$t_{\zeta}=2^{-1}[T+t_0-\rho_1].$$ Since $t_{\zeta}>t_0$, one must have $T-t_0>\rho_1$. Meanwhile, by~\eqref{type2} and our standing assumption in the current subcase, $\bar x_1(t_{\zeta})<1$. So, $T-t_0<\rho_1+2\rho_2=a^{-1}(3-x_0)$. Combining this and the inequality $T-t_0>\rho_1$ yields $\rho_1<T-t_0<a^{-1}(3-x_0)$. Our results in this subcase can be summarized as follows:
	
	$\bullet$ $\bar x_1(.)$ is given by \eqref{type1}, provided that $T-t_0=\rho_1$.
	
	$\bullet$ $\bar x_1(.)$ is given by \eqref{type2}, provided that $\rho_1<T-t_0<a^{-1}(3-x_0)$.
	
	\underline{\textit{Subcase 3b}}: ${\mathcal T}_1\neq \emptyset$. Then we have $t_0\leq \alpha_1\leq\alpha_2<T$. It follows from assertion~S2 of Proposition~\ref{Lemma_main_tool} that $T-\alpha_2=2a^{-1}$ and $\bar x_1(t)=1-a(t-\alpha_2)$ for all $t\in [\alpha_2, T]$. Thus, we have $\alpha_2=T-2a^{-1}$ and $\bar x_1(t)=1-a(t-T+2a^{-1})$ for all $t\in [T-2a^{-1}, T]$.
	 
	 If $\alpha_1<\alpha_2$, then $\bar x_1(t)=1$ for all $t\in [\alpha_1, \alpha_2]$. Indeed, suppose on the contrary that there exists $\bar t\in (\alpha_1, \alpha_2)$ satisfying $\bar x_1(\bar t)<1$. Set
	  $$\bar\alpha_1=\max\{ t\in [\alpha_1, \bar t] \;:\; \bar x_1(t)=1\} \quad {\rm and} \quad \bar\alpha_2=\min\{ t\in [\bar t, \alpha_2] \;:\; \bar x_1(t)=1\}.$$
	 Clearly, $[\bar\alpha_1, \bar\alpha_2]\subset[\alpha_1, \alpha_2]\subset [t_0, T)$ and $\bar x_1(t)<1$ for all $t\in (\bar\alpha_1,\bar\alpha_2)$. 
	  This and the condition $\bar x_1(t)>-1$ for all $t\in [t_0, T)$ imply that $\bar x_1(t)\in (-1, 1)$ for all $t\in (\bar\alpha_1,\bar\alpha_2)$. So, by assertion~S4 of Proposition~\ref{Lemma_main_tool}, we obtain a contradiction. Our claim has been proved.
	  
	If $t_0<\alpha_1$, then $\bar x_1(t)\in (-1, 1)$ for all $t\in [t_0, \alpha_1)$ and $\bar x_1(\alpha_1)=1$. Thus, repeating the arguments in the proof of assertion~S1 of Proposition~\ref{Lemma_main_tool}, we find that $\bar x_1(t)=x_0+a(t-t_0)$ for all $t\in [t_0, \alpha_1]$. As $\bar x_1(\alpha_1)$=1, we have $\alpha_1=t_0+\rho_2$. Consequently, the inequality as $T-t_0\geq(\alpha_1-t_0)+(T-\alpha_2)$ implies that $T-t_0\geq \rho_2+ 2a^{-1}=a^{-1}(3-x_0)$. Our results in this subcase can be summarized as follows:
		
		$\bullet$ $\bar x_1(.)$ is given by 
		\begin{equation*}
		\bar x_1(t)=
		\begin{cases}
		x_0+a(t-t_0), \quad & t\in [t_0, T-2a^{-1}]\\
		-1-a(t-T), & t\in (T-2a^{-1}, T],
		\end{cases}
		\end{equation*}
		provided that $T-t_0=a^{-1}(3-x_0)$.
		
		$\bullet$ $\bar x_1(.)$ is given by 
		\begin{equation*}
				\bar x_1(t)=
				\begin{cases}
				x_0+a(t-t_0), \quad & t\in [t_0, t_0+\rho_2]\\
				1, & t\in (t_0+\rho_2, T-2a^{-1}]\\
				-1-a(t-T), & t\in (T-2a^{-1}, T],
				\end{cases}
				\end{equation*}
				 provided that $T-t_0>a^{-1}(3-x_0)$.
	  
	  {\bf Case 4:}\textit{ $\bar x_1(t_0)=\bar x_1(T)=-1$ and $\bar x_1(t)>-1$ for all $t\in (t_0, T)$}.
	  
	 \underline{\textit{Subcase 4a}}: ${\mathcal T}_1=\emptyset$. Then $\bar x_1(t)\in (-1, 1)$ for all $t\in (t_0, T)$. Thus, by assertion S3 of Proposition~\ref{Lemma_main_tool} one has $T-t_0< 4a^{-1}$ and
	  \begin{equation*}
	  \bar x_1(t)=
	  \begin{cases}
	  -1+a(t-t_0), \quad & t\in [t_0, 2^{-1}(t_0+T)]\\
	  -1-a(t-T), & t\in (2^{-1}(t_0+T), T].
	  \end{cases}
	  \end{equation*}
	  
	\underline{\textit{Subcase 4b}}: ${\mathcal T}_1\neq\emptyset$. Then, the numbers $\alpha_1$ and $\alpha_2$ exist and $t_0< \alpha_1\leq\alpha_2<T$. It follows from statements S1 and S2 of Proposition~\ref{Lemma_main_tool} that $\alpha_1-t_0=T-\alpha_2=2a^{-1}$ and $\bar x_1(t)=-1+a(t-t_0)$ for all $t\in [t_0, \alpha_1]$ and $\bar x_1(t)=1-a(t-\alpha_2)$ for all $t\in [\alpha_2, T]$. Thus, we have $\alpha_1=t_0+2a^{-1}$, $\alpha_2=T-2a^{-1}$, $\bar x_1(t)=-1+a(t-t_0)$ for all $t\in [t_0, t_0+2a^{-1}]$, and $\bar x_1(t)=1-a(t-T+2a^{-1})$ for all $t\in [T-2a^{-1}, T]$. Note that one must have $T-t_0\geq 4a^{-1}$ in this subcase as $T-t_0\geq(\alpha_1-t_0)+(T-\alpha_2)$.
	
	If $T-t_0> 4a^{-1}$, i.e., $\alpha_1<\alpha_2$,  then by the result given in Subcase~3b we have $\bar x_1(t)=1$ for all $t\in [t_0+2a^{-1}, T-2a^{-1}]$.
	
	 Our results in this case can be summarized as follows:
	 	
	 	$\bullet$ $\bar x_1(.)$ is given by
	 	\begin{equation*}
	 		  \bar x_1(t)=
	 		  \begin{cases}
	 		  -1+a(t-t_0), \quad & t\in [t_0, 2^{-1}(t_0+T)]\\
	 		  -1-a(t-T), & t\in (2^{-1}(t_0+T), T],
	 		  \end{cases}
	 		  \end{equation*}
	  	provided that $T-t_0\leq 4a^{-1}$.
	  	
	 			$\bullet$ $\bar x_1(.)$ is given by 
	 		\begin{equation}
	 		\bar x_1(t)=
	 		\begin{cases}
	 		-1+a(t-t_0), \quad & t\in [t_0, t_0+2a^{-1}]\\
	 		1, & t\in (t_0+2a^{-1}, T-2a^{-1}]\\
	 		-1-a(t-T), & t\in (T-2a^{-1}, T],
	 		\end{cases}
	 		\end{equation}
	 		provided that $T-t_0>4a^{-1}$.
	 	 
	  	\medskip
	  	Now we turn our attention back to the original problem $(FP_3)$, which has a $W^{1, 1}$ global solution (see Subsection~\ref{SE_2-sided_S_constraint}). Using the given constants $a,\lambda$ with $a>\lambda>0$, we define $\rho=\dfrac{1}{\lambda}\ln \dfrac{a}{a-\lambda}>0$. This number $\rho$ is a characteristic constant of $(FP_3)$. From the analysis given in the present section we can obtain a complete synthesis of optimal processes. Due to the complexity of the possible trajectories, we prefer to present our results in six separate theorems. The first one treats the situation where $\rho\geq 2a^{-1}$, while the other five deal with the situation where $\rho< 2a^{-1}$. 
	  	
	  \begin{Theorem}\label{Thm3a}  If $\rho\geq 2a^{-1}$, then problem $(FP_3)$ has a unique local solution $(\bar x,\bar u)$, which is a unique global solution, where $\bar u(t)=-a^{-1}\dot{\bar x}(t)$ for almost everywhere $t\in [t_0, T]$ and $\bar x(t)$ can be described as follows: 
	  	  		\begin{description}
	  			\item{\rm (a)} If $T-t_0\leq a^{-1}(1+x_0)$, then 
	  			\begin{equation*}
	  				 \bar x(t)=x_0-a(t-t_0), \quad  t\in [t_0, T].
	  				 \end{equation*}
	  			\item{\rm (b)} If $a^{-1}(1+x_0)<T-t_0<a^{-1}(3-x_0)$, then 
	  			\begin{equation*}
	  				 	  \bar x(t)=
	  				 	  \begin{cases}
	  				 	  x_0+a(t-t_0), \quad & t\in [t_0, t_{\zeta}]\\
	  				 	  -1-a(t-T), & t\in (t_{\zeta}, T],
	  				 	  \end{cases}
	  				 	  \end{equation*}
	  				with $t_{\zeta}:=2^{-1}[T+t_0-a^{-1}(1+x_0)]$. 
	  			\item{\rm (c)} If $T-t_0\geq a^{-1}(3-x_0)$, then 
	  			\begin{equation*}
	  							\bar x(t)=
	  							\begin{cases}
	  							x_0+a(t-t_0), \quad & t\in [t_0, t_0+a^{-1}(1-x_0)]\\
	  							1, & t\in (t_0+a^{-1}(1-x_0), T-2a^{-1}]\\
	  							-1-a(t-T), & t\in (T-2a^{-1}, T].
	  							\end{cases}
	  							\end{equation*}
	  		\end{description}
	  	\end{Theorem}
	  	
	  	\begin{proof}  
	  	Suppose that  $\rho\geq 2a^{-1}$. Let $\rho_1, \rho_2$ be defined as in \eqref{two_rhos}. Then, one has  $\rho\geq \rho_1$, $2\rho-\rho_1\geq \rho+\rho_2$, $\rho+\rho_2\geq 2a^{-1}+\rho_2$, $2\rho \geq \rho+2a^{-1}$, and $\rho+2a^{-1}\geq 4a^{-1}$. Thus, the analysis of Case~1 and Case~2 given before this theorem tells us that the situation in Subcase~1a happens when $T-t_0<\rho_1$, while the situations in Subcase~1b, Subcase~1c, and Case~2 cannot happen. Combining the results formulated in Subcase~1a, Case~3, Case~4, and noting that the function $\bar x_1(t)$ in $(FP_{3a})$ plays the role of $\bar x(t)$ in $(FP_3)$, we obtain the assertions of the theorem.
	  	\end{proof}
	  	
	  	If $\rho < 2a^{-1}$, then the locally optimal processes of $(FP_3)$ \textit{depend greatly on the relative position of $x_0$} in the segment $[-1,1]$. In the forthcoming theorems, we distinguish five alternatives (one instance must occur, and any instance excludes others): 
	  	
	  	(i) $x_0=-1$; 
	  	
	  	(ii)  $x_0>-1$ and $\rho<a^{-1}(1+x_0)<\rho+a^{-1}(1-x_0)$; 
	  	
	  	(iii) $x_0>-1$ and $a^{-1}(1+x_0)=\rho+a^{-1}(1-x_0)$;
	  	
	  	(iv) $x_0>-1$ and $a^{-1}(1+x_0)>\rho+a^{-1}(1-x_0)$;
	  	
	  	(v) $x_0>-1$ and $a^{-1}(1+x_0)\leq\rho$.
	  	
	  	\medskip
	  	It is worthy to stress that to describe the possibilities (i)--(v) we have used just the parameters $a,\lambda,$ and $x_0$. In each one of the situations (i)--(v), the synthesis of the trajectories suspected for local minimizers of $(FP_3)$ is obtained by considering \textit{the position of the number $T-t_0>0$ on the half-line $[0,+\infty)$}, which is divided into sections by the values $\rho$, $2\rho$, $\rho+2a^{-1}$, $4a^{-1}$, and other constants appeared in (i)--(v).
	  	
	   	\begin{Theorem}\label{Thm3b}  If $\rho < 2a^{-1}$ and $x_0=-1$, then any local solution of problem $(FP_3)$ must have the form $(\bar x,\bar u)$, where $\bar u(t)=-a^{-1}\dot{\bar x}(t)$ for almost everywhere $t\in [t_0, T]$ and $\bar x(t)$ is described as follows: 
	  		\begin{description}
	  		\item{\rm (a)} If $T-t_0\leq 2\rho$, then 
	  		\begin{equation}\label{Thm3b-traj1}
	  			 		  \bar x(t)=
	  			 		  \begin{cases}
	  			 		  -1+a(t-t_0), \quad & t\in [t_0, 2^{-1}(t_0+T)]\\
	  			 		  -1-a(t-T), & t\in (2^{-1}(t_0+T), T].
	  			 		  \end{cases}
	  			 		  \end{equation}
	  			\item{\rm (b)} If $2\rho< T-t_0<\rho+2a^{-1}$, then $\bar x(t)$ is given by either \eqref{Thm3b-traj1}, or
	  			\begin{equation*}
	  				  	 	  	    		 \bar x(t)=
	  				  	 	  	    		 \begin{cases}
	  				  	 	  	    		 -1+a(t-t_0), \quad & t \in [t_0, \bar t]\\
	  				  	 	  	    		 -1-a(t+t_0-2\bar t), & t \in (\bar t, T],
	  				  	 	  	    		 \end{cases}
	  				  	 	  	    		 \end{equation*} where $\bar t=T-\rho$. 
	  			\item{\rm (c)} If $\rho+2a^{-1}\leq T-t_0<4a^{-1}$, then $\bar x(t)$ is given  by either \eqref{Thm3b-traj1}, or
	  			\begin{equation*}
	  			  	  	    		 \bar x(t)=
	  			  	  	    		 \begin{cases}
	  			  	  	    		  -1+a(t-t_0), \quad & t \in [t_0, t_0+2a^{-1}]\\
	  			  	  	    		 1-a(t-t_0-2a^{-1}), & t \in (t_0+2a^{-1}, T].
	  			  	  	    		 \end{cases}
	  			  	  	    		 \end{equation*}
	  			  	  		\item{\rm (d)}	If $T-t_0=4a^{-1}$, then $\bar x(t)$ is given by \eqref{Thm3b-traj1}. 
	  			  		\item{\rm (e)}	If $T-t_0>4a^{-1}$, then 
	  				  			  		    \begin{equation*}
	  				  			  		    \bar x(t)=
	  				  			  			\begin{cases}
	  				  			  			-1+a(t-t_0), \quad & t\in [t_0, t_0+2a^{-1}]\\
	  				  			  		    1, & t\in (t_0+2a^{-1}, T-2a^{-1}]\\
	  				  			  		    -1-a(t-T), & t\in (T-2a^{-1}, T].
	  				  			  			\end{cases}
	  				  			  			\end{equation*}
	  				  		\end{description}
	  				  		
	  		In the situations described in $\rm (a), \rm (d)$, and $\rm (e)$, $(\bar x,\bar u)$ is a unique local solution  of $(FP_3)$, which is also a unique global solution of the problem.
	  	\end{Theorem}
	  	\begin{proof} Suppose that  $\rho< 2a^{-1}$ and $x_0=-1$. To obtain the assertions (a)--(e), it suffices to combine the results formulated in Case~2 and Case~4, having in mind that $\bar x_1(t)$ in $(FP_{3a})$ plays the role of $\bar x(t)$ in $(FP_3)$.
	  	\end{proof}
	  	
	 \begin{Theorem}\label{Thm3c1}  If $\rho < 2a^{-1}$, $x_0>-1$, and $\rho<a^{-1}(1+x_0)<\rho+a^{-1}(1-x_0)$, then any local solution of problem $(FP_3)$ must have the form $(\bar x,\bar u)$, where $\bar u(t)=-a^{-1}\dot{\bar x}(t)$ for almost everywhere $t\in [t_0, T]$ and $\bar x(t)$ is described as follows: 
	  		  		\begin{description}
	  		  		\item{\rm (a)} If $T-t_0\leq \rho$, then $\bar x(t)$ is given by
	  		  		\begin{equation}\label{Thm3c-traj1}
	  		  	    \bar x(t)=x_0-a(t-t_0), \quad  t\in [t_0, T].
	  		  		\end{equation}
	  		         \item{\rm (b)}	If $\rho<T-t_0<a^{-1}(1+x_0)$, then 
	  		         \begin{equation}\label{Thm3c-traj2}
	  		         	 \bar x(t)=
	  		         	 \begin{cases}
	  		         	 x_0+a(t-t_0), \quad & t \in [t_0, \bar t]\\
	  		         	 x_0-a(t+t_0-2\bar t), & t \in (\bar t, T],
	  		         	 \end{cases}
	  		         \end{equation} where $\bar t=T-\rho$. 
	  		        \item{\rm (c)} 	If $T-t_0=a^{-1}(1+x_0)$, then  $\bar x(t)$ is given by either \eqref{Thm3c-traj1}, or \eqref{Thm3c-traj2}.
	  		         \item{\rm (d)} 	If $a^{-1}(1+x_0)<T-t_0<\rho+a^{-1}(1-x_0)$, then $\bar x(t)$ is given by either \eqref{Thm3c-traj2}, or 
	  		  	     \begin{equation}\label{Thm3c-traj4}
	  		  		 \bar x(t)=
	  		  		 \begin{cases}
	  		  		 x_0+a(t-t_0), \quad & t\in [t_0, t_{\zeta}]\\
	  		  		 -1-a(t-T), & t\in (t_{\zeta}, T],
	  		  		 \end{cases}
	  		  		 \end{equation}
	  		  		  with $t_{\zeta}:=2^{-1}[T+t_0-a^{-1}(1+x_0)]$. 		
	  		  	  \item{\rm (e)} If $\rho+a^{-1}(1-x_0)\leq T-t_0<a^{-1}(3-x_0)$, then $\bar x(t)$ is given by either \eqref{Thm3c-traj4}, or
	  		  	 \begin{equation}\label{Thm3c-traj3}
	  		  	\bar x(t)=
	  		  	\begin{cases}
	  		  	x_0+a(t-t_0), \quad & t \in [t_0, t_0+a^{-1}(1-x_0)]\\
	  		  	1-a(t-t_0-a^{-1}(1-x_0)), & t \in (t_0+a^{-1}(1-x_0), T].
	  		  	\end{cases}
	  		    \end{equation}
	  	\item{\rm (f)} If $T-t_0 = a^{-1}(3-x_0)$, then
	  		 \begin{equation}\label{Thm3c-traj4a}
	  		 \bar x(t)=
	  		 \begin{cases}
	  		 x_0+a(t-t_0), \quad & t\in [t_0, T-2a^{-1}]\\
	  		 -1-a(t-T), & t\in (T-2a^{-1}, T].
	  		 \end{cases}
	  		 \end{equation}
	  	  	\item{\rm (g)} If $T-t_0> a^{-1}(3-x_0)$, then 
	  	\begin{equation}\label{Thm3c-traj5}
	     \bar x(t)=
	     \begin{cases}
	  	x_0+a(t-t_0), \quad & t\in [t_0, t_0+a^{-1}(1-x_0)]\\
	  	1, & t\in (t_0+a^{-1}(1-x_0), T-2a^{-1}]\\
	  	-1-a(t-T), & t\in (T-2a^{-1}, T].
	   \end{cases}
	   \end{equation}
	  		 \end{description}
	  		 
	  		 	In the situations described in $\rm (a), \rm (b), \rm (f)$, and $\rm (g)$, $(\bar x,\bar u)$ is a unique local solution of~$(FP_3)$, which is also a unique global solution of the problem.
	  		 \end{Theorem}
	  		  \begin{proof}
	  		 Suppose that  $\rho < 2a^{-1}$, $x_0>-1$, and $\rho<a^{-1}(1+x_0)<\rho+a^{-1}(1-x_0)$. Let $\rho_1, \rho_2$ be defined as in \eqref{two_rhos}. Then, combining the results formulated in Case~1 and Case~3, and noting that the function $\bar x_1(t)$ in $(FP_{3a})$ plays the role of $\bar x(t)$ in $(FP_3)$, we obtain the assertions (a) -- (g).
	  		 \end{proof}
	  		 
\begin{Theorem}\label{Thm3c2} If $\rho < 2a^{-1}$, $x_0>-1$, and $a^{-1}(1+x_0)=\rho+a^{-1}(1-x_0)$, then any local solution of problem $(FP_3)$ must have the form $(\bar x,\bar u)$, where $\bar u(t)=-a^{-1}\dot{\bar x}(t)$ for almost everywhere $t\in [t_0, T]$ and $\bar x(t)$ is described as follows: 
\begin{description}
\item{\rm (a)} If $T-t_0\leq \rho$, then $\bar x(t)$ is given by \eqref{Thm3c-traj1}. 
\item{\rm (b)}	If $\rho<T-t_0<a^{-1}(1+x_0)$, then $\bar x(t)$ is given by \eqref{Thm3c-traj2}. 
\item{\rm (c)} 	If $T-t_0=a^{-1}(1+x_0)$, then $\bar x(t)$ is given by either  \eqref{Thm3c-traj1}, or \eqref{Thm3c-traj3}.
\item{\rm (d)} 	If $a^{-1}(1+x_0)<T-t_0<a^{-1}(3-x_0)$, then $\bar x(t)$ is given by either \eqref{Thm3c-traj4}, or \eqref{Thm3c-traj3}.
\item{\rm (f)} If $T-t_0 = a^{-1}(3-x_0)$, then  $\bar x(t)$ is given by \eqref{Thm3c-traj4a}. 
\item{\rm (g)} If $T-t_0> a^{-1}(3-x_0)$, then $\bar x(t)$ is given by \eqref{Thm3c-traj5}. 
\end{description}

In the situations described in $\rm (a), \rm (b), (f)$, and $\rm (g)$, $(\bar x,\bar u)$ is a unique local solution  of~$(FP_3)$, which is also a unique global solution of the problem.
\end{Theorem}
\begin{proof}
Suppose that $\rho < 2a^{-1}$, $x_0>-1$, and $a^{-1}(1+x_0)=\rho+a^{-1}(1-x_0)$. Then, combining the results formulated in Case~1 and Case~3, and noting that the function $\bar x_1(t)$ in $(FP_{3a})$ plays the role of $\bar x(t)$ in $(FP_3)$, we obtain the desired assertions.
\end{proof}	
 
\begin{Theorem}\label{Thm3c3}  If $\rho < 2a^{-1}$, $x_0>-1$, and $a^{-1}(1+x_0)>\rho+a^{-1}(1-x_0)$, then any local solution of problem $(FP_3)$ must have the form $(\bar x,\bar u)$, where $\bar u(t)=-a^{-1}\dot{\bar x}(t)$ for almost everywhere $t\in [t_0, T]$ and $\bar x(t)$ is described as follows: 
\begin{description}
\item{\rm (a)} If $T-t_0\leq \rho$, then $\bar x(t)$ is given by \eqref{Thm3c-traj1}. 
\item{\rm (b)}	If $\rho<T-t_0<\rho+a^{-1}(1-x_0)$, then $\bar x(t)$ is given by \eqref{Thm3c-traj2}. 
\item{\rm (c)} 	If $\rho+a^{-1}(1-x_0)<T-t_0<a^{-1}(1+x_0)$, then $\bar x(t)$ is given by \eqref{Thm3c-traj3}. 
\item{\rm (d)} 	If $T-t_0=a^{-1}(1+x_0)$, then $\bar x(t)$ is given by either  \eqref{Thm3c-traj1}, or \eqref{Thm3c-traj3}.
\item{\rm (e)} 	If $a^{-1}(1+x_0)<T-t_0<a^{-1}(3-x_0)$, then $\bar x(t)$ is given by either \eqref{Thm3c-traj4}, or \eqref{Thm3c-traj3}.
\item{\rm (f)} If $T-t_0 = a^{-1}(3-x_0)$, then $\bar x(t)$ is given by \eqref{Thm3c-traj4a}. 
\item{\rm (g)} If $T-t_0> a^{-1}(3-x_0)$, then $\bar x(t)$ is given by \eqref{Thm3c-traj5}. 
\end{description}

In the situations described in $\rm (a), \rm (b), \rm (c), \rm(f)$, and $\rm (g)$, $(\bar x,\bar u)$ is a unique local solution  of~$(FP_3)$, which is also a unique global solution of the problem.
\end{Theorem}
\begin{proof}
Suppose that $\rho < 2a^{-1}$, $x_0>-1$, and $a^{-1}(1+x_0)>\rho+a^{-1}(1-x_0)$. Then, combining the results formulated in Case~1 and Case~3, and noting that the function $\bar x_1(t)$ in $(FP_{3a})$ plays the role of $\bar x(t)$ in $(FP_3)$, we obtain the assertions of the theorem.
\end{proof}
 	  		 	  		 	  		  	
\begin{Theorem}\label{Thm3d}  If $\rho < 2a^{-1}$, $x_0>-1$, and $a^{-1}(1+x_0)\leq\rho$, then any local solution of problem $(FP_3)$ must have the form $(\bar x,\bar u)$, where $\bar u(t)=-a^{-1}\dot{\bar x}(t)$ for almost everywhere $t\in [t_0, T]$ and $\bar x(t)$ can be described as follows: 
\begin{description}
\item{\rm (a)} If $T-t_0\leq a^{-1}(1+x_0)$, then $\bar x(t)$ is given by \eqref{Thm3c-traj1}. 
\item{\rm (b)} 	If $a^{-1}(1+x_0)<T-t_0<2\rho-a^{-1}(1+x_0)$, then $\bar x(t)$ is given by \eqref{Thm3c-traj4}. 
\item{\rm (c)} 	If $2\rho-a^{-1}(1+x_0)<T-t_0<\rho+a^{-1}(1-x_0)$, then $\bar x(t)$ is given by either \eqref{Thm3c-traj2}, or \eqref{Thm3c-traj4}.
\item{\rm (d)} If $\rho+a^{-1}(1-x_0)\leq T-t_0<a^{-1}(3-x_0)$, then $\bar x(t)$ is given by either \eqref{Thm3c-traj4}, or \eqref{Thm3c-traj3}.
\item{\rm (e)} If $T-t_0 = a^{-1}(3-x_0)$, then $\bar x(t)$ is given by  \eqref{Thm3c-traj4a}. 
\item{\rm (f)} If $T-t_0> a^{-1}(3-x_0)$, then $\bar x(t)$ is given by \eqref{Thm3c-traj5}. 
\end{description}

In the situations described in $\rm (a), \rm (b), \rm(e)$, and $\rm (f)$, $(\bar x,\bar u)$ is a unique local solution  of~$(FP_3)$, which is also a unique global solution of the problem.
\end{Theorem}
\begin{proof}
Suppose that $\rho < 2a^{-1}$, $x_0>-1$, and $a^{-1}(1+x_0)\leq\rho$. Let $\rho_1, \rho_2$ be given by \eqref{two_rhos}. Then, combining the results formulated in Case~1 and Case~3, and noting that the function $\bar x_1(t)$ in $(FP_{3a})$ plays the role of $\bar x(t)$ in $(FP_3)$, we obtain the assertions (a) -- (f).
\end{proof}
	  	
	    \section{Conclusions}\label{Conclusions}
	    
	    We have analyzed a maximum principle for finite horizon state constrained problems via one parametric example of optimal control problems of the Langrange type, which has five parameters. This problem resembles the optimal growth problem in mathematical economics. It belongs to the class of control problems with bilateral state constraints. We have proved that the control problem in the example can have not more than  two local solutions, and at least one of them  which must be a global solution. Moreover, we have presented explicit descriptions of the optimal processes, which are suspected to be local solutions, with respect to the five parameters.
	    
	    The obtained results allow us to have a deep understanding of the maximum principle in question. 
	   	    
	   	It seems to us that economic optimal growth models can be studied by advanced tools from functional analysis and optimal control theory via the approach adopted in this paper.

\end{document}